\documentclass[12pt]{amsart}
\usepackage{amsfonts, amsmath, amstext, latexsym,amssymb}
\usepackage[english]{babel}
\usepackage{graphpap}  
\usepackage[all]{xy}

\begin{document}  

\newcommand{\rr}{\ensuremath{\mathbf R}}

\newcommand{\cccc}{\ensuremath{\mathcal{C}}}
\newcommand{\ff}{\ensuremath{\mathcal{F}}}
\newcommand{\uu}{\ensuremath{\mathcal{U}}}
\newcommand{\G}{\ensuremath{\mathcal{G}}}
\newcommand{\GG}{\ensuremath{\overline{\mathcal{G}}}}

\newcommand{\al}{\ensuremath{\alpha}}
\newcommand{\bb}{\ensuremath{\beta}}
\newcommand{\be}{\ensuremath{\beta}}
\newcommand{\g}{\ensuremath{\gamma}}
\newcommand{\de}{\ensuremath{\delta}}
\newcommand{\D}{\ensuremath{\Delta}}
\newcommand{\f}{\ensuremath{\varphi}}
\newcommand{\s}{\ensuremath{\sigma}}  
\newcommand{\la}{\ensuremath{\lambda}}
\newcommand{\w}{\ensuremath{\omega}}  
\newcommand{\W}{\ensuremath{\Omega}}
\newcommand{\e}{\ensuremath{\epsilon}}
\newcommand{\ep}{\ensuremath{\varepsilon}}

\newcommand{\p}{\ensuremath{\partial}}
\newcommand{\pp}{\ensuremath{\parallel}} 
\newcommand{\inv}{\ensuremath{^{-1}}} 
\newcommand{\dhat}{\ensuremath{\hat{d}}}
\newcommand{\dbar}{\ensuremath{\bar{d}}}
\newcommand{\dcheck}{\ensuremath{\check{d}}}
\newcommand{\se}{\ensuremath{\subseteq}}
\newcommand{\lra}{\ensuremath{\longrightarrow}}
\newcommand{\lla}{\ensuremath{\longleftarrow}}
\newcommand{\dd}{\ensuremath{\hat{d}}}
\newcommand{\ds}{\ensuremath{\displaystyle}}
\newcommand{\CAT}{\ensuremath{\mbox{CAT}}}
\newcommand{\lbr}{\ensuremath{ | \! [ }}
\newcommand{\rbr}{\ensuremath{ ] \! | }} 

\newcommand{\go}{\ensuremath{{\mathcal{G}}^{(0)}}}
\newcommand{\fff}{\ensuremath{\bar{f}}}
\newcommand{\?}{{\bf????????}}

\def\weak{weak-* }
\def\ti{-\allowhyphens}
\def\hb{{\rm H}_{\rm b}}
\def\hc{{\rm H}_{\rm c}}
\def\hbc{{\rm H}_{\rm cb}}
\def\coo{C_{00}}
\def\co{C_0}
\def\No{N\raise4pt\hbox{\tiny o}\kern+.2em}
\def\ro{\varrho}
\def\isom{\mathrm{Isom}}
\def\fhi{\varphi}

\newtheoremstyle{mydefinition}% name
  {}%      Space above, empty = `usual value'
  {}%      Space below
  {\upshape}% Body font
  {}%         Indent amount (empty = no indent, \parindent = para indent)
  {\bf}% Thm head font
  {.}%        Punctuation after thm head
  { }%     Space after thm head: " " = normal interword space;
        %       \newline = linebreak
  {}% Thm head spec

\newtheorem{ttt}{Theorem}  
\newtheorem*{ttt*}{Theorem}
\newtheorem{lem}[ttt]{Lemma}
\newtheorem{slem}[ttt]{Sublemma}
\newtheorem{ddd}[ttt]{Definition}  
\newtheorem{ppp}[ttt]{Proposition}  
\newtheorem{ccc}[ttt]{Corollary}
\newtheorem*{ccc*}{Corollary}
\theoremstyle{mydefinition}
\newtheorem{rem}[ttt]{Remark} 

\newtheorem{ttt-in-text}{Theorem}
%\renewcommand{\baselinestretch}{1.3} 
%\small\normalsize  

\title[ideal bicombings]{Ideal bicombings for hyperbolic groups\\ and applications}
\author{Igor Mineyev, Nicolas Monod and Yehuda Shalom}
\thanks{I.M. is~partially supported by NSF grant DMS~0132514 (CAREER);
        N.M. is~partially supported by FNS grant~8220-067641 and NSF grant DMS~0204601.}
%\date{}
\begin{abstract}
For every hyperbolic group and more general hyperbolic graphs, we construct an equivariant ideal bicombing: this is a homological analogue of the geodesic flow on negatively curved manifolds. We then construct a cohomological invariant which implies that several Measure Equivalence and Orbit Equivalence rigidity results established in~\cite{Monod-Shalom2} hold for all non-elementary hyperbolic groups and their non-elementary subgroups. We also derive superrigidity results for actions of general irreducible lattices on a large class of hyperbolic metric spaces.
\end{abstract}
\maketitle 
%\renewcommand{\thefootnote}{\fnsymbol{footnote}}
%\footnote[0]{\textit{Date}: }
%\footnote[0]{1991 \textit{Mathematics Subject Classification}:
%}

%%%%%%%%%%%%%%%%%%%%%%%%%%%%%%%%%%%%%%%%%%%%%%%%%%%%%%%%%%%%%%%%%%
\section{Introduction}

\noindent
\textbf{1.1.} If $M$ is a simply connected manifold with sectional curvature $\leq -1$, then any two distinct points of the sphere at infinity $\partial M$ can be connected by a unique geodesic line. Moreover, when two such geodesics share a common ideal endpoint, they converge at exponential synchronous rate, and the whole setting is of course equivariant under the group of isometries of $M$. A similar situation occurs more generally for $\CAT(-1)$ spaces. \itshape The goal of this paper is to establish an analogous phenomenon for Gromov-hyperbolic graphs, and to use it in constructing cohomological invariants that have far-reaching applications in rigidity theory.\upshape

One of the reasons the theory of hyperbolic groups is so fruitful is that whilst hyperbolicity is a very general and robust property (indeed, it is \emph{generic} in some ways~\cite{Gromov91,Olshanskii92,ZukERN}), many ideas from large scale differential geometry in negative curvature do find analogues for hyperbolic groups. In the case at hand, geodesics are to be replaced with a homological concept, an \emph{ideal} version of the notion of homological bicombing:

Let $\G=(\go,E)$ a (connected) hyperbolic graph, for instance a Cayley graph of a hyperbolic group, and let $\p\G$ be its boundary at infinity (for all terminology see Section~\ref{sec_notation}).

\begin{ddd}
\label{def_ideal}%
An {\sf ideal bicombing} on $\G$ is a map $q:(\p\G)^2\to \ell^\infty(E)$ such that:

\begin{itemize}
\item[(i)] $q[\xi,\eta]$ is a cycle for all $\xi,\eta\in\p\G$, \emph{i.e.} $\partial q[\xi,\eta]=0$.
\item[(ii)] For all distinct $\xi_+,\xi_-\in\p\G$ and all disjoint graph neighbourhoods $V_\pm$ of $\xi_\pm$, the function $\partial (q[\xi_-, \xi_+]|_{V_\pm})$ is finitely supported of sum $\mp 1$.
\end{itemize}
\end{ddd}

We say that $q$ has {\sf bounded area} if $|\!|q[\xi,\eta]+q[\eta,\zeta]+q[\zeta,\xi]{|\!|}_1$ is finite and uniformly bounded over $\xi,\eta,\zeta\in\p\G$, and that $q$ is {\sf quasigeodesic} if there exists a constant $C$ such that for all $\xi,\eta\in\p\G$ the support of~$q[\xi,\eta]$ lies in the $C$\ti neighbourhood of any geodesic from $\xi$ to $\eta$.

\medskip

Our first result is obtained by proving that (a generalization of) the homological bicombing constructed in~\cite{Mineyev01} can be continuously extended to infinity; we show:

\begin{ttt}
\label{thm_existence}%
Every hyperbolic graph $\G$ of bounded valency admits a \weak continuous quasi-geodesic ideal bicombing of bounded area which is equivariant under the automorphism group of $\G$.
\end{ttt}

This result is our replacement for infinite geodesics with synchronous convergence. We point out that Gromov introduced a geodesic flow for hyperbolic groups (see~\cite{Gromov87,Champetier}), but his construction is not sufficiently well understood for our further needs (even in that case). We single out the precise notion of (homological) exponential convergence in Theorem~\ref{aa'bb'} below.

\medskip

\noindent
\textbf{1.2.} Our main motivation for Theorem~\ref{thm_existence} is the method introduced in~\cite{Monod-Shalom1}, where an invariant in bounded cohomology is associated to the boundary at infinity of any proper $\CAT(-1)$ space. Exploiting the theory developed in~\cite{Burger-Monod3,Monod}, this invariant is used in~\cite{Monod-Shalom1} to prove superrigidity results and applied in~\cite{Monod-Shalom2} to the theory of Orbit Equivalence to establish new rigidity phenomena. Of particular importance in that context is the class~$\mathcal{C}_\mathrm{reg}$ of groups $\Gamma$ for which $\hb^2(\Gamma, \ell^2(\Gamma))$ is non-trivial; it is shown in~\cite{Monod-Shalom1} that this class contains all groups acting properly and non-elementarily on proper $\CAT(-1)$ spaces, various amalgamated products (including those with finite amalgamated subgroup), and hyperbolic groups with vanishing first $\ell^2$\ti Betti number. We can now complete the picture as conjectured in~\cite{Monod-Shalom1}:

\begin{ttt}
\label{thm_coho}%
Let $\Gamma$ be a countable group admitting a proper non-elementary action on a hyperbolic graph of bounded valency. Then $\hb^2(\Gamma, \ell^p(\Gamma))$ is non-trivial for all $1\leq p <\infty$. In particular $\Gamma$ is in the class~$\mathcal{C}_\mathrm{reg}$.
\end{ttt}

This applies in particular to the case where $\Gamma$ is a non-elementary hyperbolic group or more generally a non-elementary subgroup of a hyperbolic group. We point out that whilst there are countably many non-isomorphic hyperbolic groups, there is a continuum~$2^{\aleph_0}$ of non-isomorphic subgroups of hyperbolic groups (a fact explained to us by I.~Kapovich and P.~Schupp, see Remark~\ref{rem_Ilya}); notice that only countably many of them are elementary. This adds a substantial collection of groups to the class~$\mathcal{C}_\mathrm{reg}$ (which was already shown to have cardinality~$2^{\aleph_0}$ in~\cite{Monod-Shalom1} since it contains the free product of any two countable groups).

\smallskip

Using structure theory of locally compact groups, one can combine Theorem~\ref{thm_coho} with the case of $\CAT(-1)$ spaces addressed in~\cite{Monod-Shalom1} and deduce the following:

\begin{ccc}
\label{cor_general_hyp}%
Let $X$ be a hyperbolic proper geodesic metric space on which $\isom(X)$ acts cocompactly. Then any countable group $\Gamma$ admitting a proper non-elementary isometric action on $X$ is in~$\mathcal{C}_\mathrm{reg}$.
\end{ccc}

Theorem~\ref{thm_coho} is obtained by first constructing a specific cocycle at infinity
$$\omega:\ (\p\G)^3 \longrightarrow \ell^1(E\times E)$$
\noindent
using a ``doubling'' of the ideal bicombing, see Theorem~\ref{thm_cocycle} below. This cocycle is our cohomological invariant mentioned above. Then we obtain a class in bounded cohomology using the functorial approach~\cite{Burger-Monod3,Monod} together with Poisson boundary theory. These elements allow us to argue in analogy with~\cite{Monod-Shalom1}.

\medskip

\noindent
\textbf{1.3.} Our interest in Theorem~\ref{thm_coho} lies firstly in the attempt to relate more intimately geometric and measurable group theory (as \emph{e.g.} in~\cite{GaboriauL2,Monod-Shalom2,ShalomQI}). Recall that two countable groups $\Gamma$ and $\Lambda$ are called {\sf measure equivalent}~(ME) if there is a $\sigma$\ti finite measure space $(\Sigma,m)$ with commuting measure preserving $\Gamma$\ti\ and $\Lambda$\ti actions, such that each one of the actions admits a finite measure fundamental domain. For example, it is evident from the definition that any two lattices in the same locally compact group are ME. The ME relation was suggested by Gromov~\cite[0.5.E]{Gromov91} as a measurable analogue of being quasi-isometric: Indeed, the latter is equivalent to the existence of a locally compact space with commuting continuous $\Gamma$\ti\ and $\Lambda$\ti actions, such that each one is proper and cocompact~\cite{Gromov91}. Thus ME amounts to replacing such a space by the measurable counterpart $(\Sigma,m)$.

\smallskip

Whilst the question of which group properties are geometric (\emph{i.e.} preserved under the quasi-isometry relation) is a well studied one, and of central interest in geometric group theory~--- much less is known about it in the analogous measurable setting, which seems to require  more subtle techniques in general. An intriguing example is provided  by the distinguished class of hyperbolic groups. Indeed, whilst being hyperbolic is well known to be a geometric property, it is \emph{not} ME invariant. For instance, in the automorphism group of a regular tree one can find non-Abelian free groups as cocompact lattices, but also non-finitely generated groups as non-uniform lattices. Hence, being ME to a hyperbolic (or even free) group does not even guarantee finite generation (let alone finite presentation or hyperbolicity). Tamer counter-examples arise when considering a non-uniform lattice in a rank one simple Lie group $G$, which is ME to any uniform (hence hyperbolic) lattice in the same $G$, but is never hyperbolic itself when $G\neq \mathrm{SL}_2(\bf{R})$. Thus, it is natural to inquire what it is in the geometric hyperbolicity property which can be, after all, detected measure-theoretically. Previously it was shown by Adams that a group ME to a hyperbolic group must have a finite centre~\cite{Adams95}, and that it is not isomorphic to a product of infinite groups~\cite{Adams94a} (see also Zimmer's~\cite{Zimmer83}). In~\cite{Monod-Shalom2} it was shown that the property of being in the class~$\mathcal{C}_\mathrm{reg}$ is a ME invariant. This, together with Theorem~\ref{thm_coho}, yields the following consequence:

\begin{ccc}
If $\Gamma$ is a countable group which is ME to a non-elementary (subgroup of a) hyperbolic group, then $\Gamma$ is in the class~$\mathcal{C}_\mathrm{reg}$.\hfill\qedsymbol
\end{ccc}

This implies in particular the two results of Adams mentioned above: If, for the first case, $\Gamma$ has a normal amenable subgroup $N$ (\emph{e.g.} its centre), then \emph{any} dual $\Gamma$\ti module with non-vanishing $\hb^2$ must have a non-zero vector invariant under $N$; see \emph{e.g.}~\cite{Noskov} or~\cite[8.5.4]{Monod}. In the case $\Gamma=\prod\Gamma_i$, any separable dual $\Gamma$\ti module with non-vanishing $\hb^2$ must have a non-zero vector invariant under some $\Gamma_i$; see \emph{e.g.}~\cite[Thm.~14]{Burger-Monod3} and~\cite{Kaimanovich03}. This prevents $\Gamma$ from being in~$\mathcal{C}_\mathrm{reg}$ once $N$, respectively $\Gamma_i$, is infinite.

In fact, these arguments give a stronger algebraic restriction on countable groups $\Gamma$ which are ME to a non-elementary (subgroup of a) hyperbolic group: Any infinite normal subgroup of $\Gamma$ has finite amenable radical and cannot split as direct product of two infinite groups. Indeed, it is shown in~\cite{Monod-Shalom2} that the class~$\mathcal{C}_\mathrm{reg}$ is stable under passing to infinite normal subgroups. (The normality assumption is essential as can be readily seen on the example of ${\bf Z}*{\bf Z}^2$, which is ME to the hyperbolic group ${\bf Z}*{\bf Z}$.)

\medskip

The following application of Theorem~\ref{thm_coho}, for which we again refer to~\cite{Monod-Shalom2}, shows a ``unique factorization'' phenomenon for torsion-free hyperbolic groups with respect to ME:

\begin{ccc}
\label{cor_prime}%
Let $\Gamma _i$ and $\Lambda_j$ be non-elementary torsion-free hyperbolic groups ($1 \le i \le n$, $1 \le j \le m$). If $\Gamma = \prod \Gamma_i$ is ME to $\Lambda =\prod \Lambda_j$ then $n=m$ and after permutation of the indices each $\Gamma_i$ is ME to $\Lambda_i$.\hfill\qedsymbol
\end{ccc}

The \emph{geometric} counterpart of this statement (prime factorization for quasi-isometries instead of ME) follows from arguments of~\cite{Eskin-Farb1,Eskin-Farb2,Kleiner-Leeb}. We emphasize that, outside the class~$\mathcal{C}_\mathrm{reg}$, the unique factorization of Corollary~\ref{cor_prime} can fail both because $n\neq m$ can occur and also when $n=m$ is assumed; see~\cite{Monod-Shalom2}.

\medskip

\noindent
\textbf{1.4.} As further consequence, the collection of measurable Orbit Equivalence rigidity results established in~\cite{Monod-Shalom2} for products of groups in the class~$\mathcal{C}_\mathrm{reg}$ now applies, by Theorem~\ref{thm_coho}, to hyperbolic groups aswell. Indeed the OE results follow in fact from the ME results. We just indicate for the record a very particular case which is easy to state:

\begin{ttt*}[\cite{Monod-Shalom2}]
Let $\Gamma = \Gamma_1 \times \Gamma_2$ where each $\Gamma_i$ is a countable torsion-free group in the class~$\mathcal{C}_\mathrm{reg}$. If a Bernoulli $\Gamma$\ti action is orbit equivalent to a Bernoulli $\Lambda$\ti action for some {\bf arbitrary} countable group $\Lambda$, then $\Gamma$ and $\Lambda$ are isomorphic.

\nobreak
Moreover, with respect to some isomorphism $\Gamma\cong\Lambda$ the actions are isomorphic via a Borel isomorphism which induces the given orbit equivalence.\hfill\qedsymbol
\end{ttt*}

See also~\cite{Kechris-Hjorth} for related results in a somewhat different setting.

\medskip

\noindent
\textbf{1.5.} Recall that a lattice $\Lambda<G=G_1\times \cdots \times G_n$ in a product of ($n\geq 2$) locally compact groups $G_i$ is {\sf irreducible} if its projection to each factor $G_i$ is dense. This generalizes the usual notion of irreducibility for lattices in semisimple algebraic groups, where each $G_i$ is simple. In this algebraic context, where $\Lambda$ must be \emph{arithmetic} by a fundamental theorem of Margulis, it is known that every homomorphism $\fhi:\Lambda\to\Gamma$ to a hyperbolic group $\Gamma$ has finite image. More generally in fact, if each $G_i$ has finite dimensional $\hbc^2(G_i,{\bf R})$, then every homomorphism to a hyperbolic group is {\sf elementary}, \emph{i.e.} has elementary image. This follows from combining the vanishing theorem of~\cite{Burger-Monod3} with the non-vanishing of~\cite{Fujiwara98}, and holds for more general hyperbolic graphs (as those in~\cite{Bestvina-Fujiwara}).

\smallskip

However, even \emph{without any assumption} on the factors of the ambient group, we now show that non-elementary actions extend after possibly factoring out a compact obstruction. This applies to situations where there are in fact non-elementary actions. 

\begin{ttt}
\label{thm_superrigidity}%
Let $G_1,\ldots, G_n$ be locally compact $\sigma$\ti compact groups, $\Lambda$ an irreducible lattice in $G=G_1\times \cdots \times G_n$ and $\pi:\Lambda\to\isom(X)$ a homomorphism, where $X$ is either:

\begin{itemize}
\item[(i)] a hyperbolic graph of bounded valency, or

\item[(ii)] any hyperbolic proper geodesic metric space on which $\isom(X)$ acts cocompactly.
\end{itemize}

\nobreak
Then either the closure $H=\overline{\pi (\Lambda)}<\isom(X)$ is amenable (and hence the $\Lambda$\ti action on $X$ is elementary), or there is a compact normal subgroup $K\lhd H$ such that the induced homomorphism $\Lambda\to H/K$ extends to a continuous homomorphism $G\to H/K$, factoring through some $G_i$.
\end{ttt}

To put the result in a better perspective, we note the following comments:

\smallskip\noindent
\textbf{a)} Observe that Theorem~\ref{thm_superrigidity} applies immediately to homomorphisms to hyperbolic groups. In particular, \itshape any non-elementary homomorphism from $\Lambda$ to a torsion-free hyperbolic group extends continuously to $G$.\upshape

\smallskip\noindent
\textbf{b)} Note that if $\Lambda$ (or equivalently each $G_i$) has in addition Kazhdan's property~(T), then it follows that without any further assumption $\pi$ {\bf always} extends continuously (modulo a compact subgroup).

\smallskip\noindent
\textbf{c)} It is also possible to derive a cocycle superrigidity result \`a la Zimmer from the same cohomological information we obtain here, which actually implies the theorem above, both for product of general groups and for (lattices in) higher rank simple Lie/algebraic groups. We refer the reader to~\cite{Monod-Shalom1} where this is explained, along with other general results.

\medskip

A worth noting feature of Theorem~\ref{thm_superrigidity} is that it goes beyond the setting of $\CAT(-1)$ target spaces, to the more robust framework of hyperbolicity. For example, in~\cite{Monod-Shalom1} it is observed that if $\Gamma$ is any Gromov-hyperbolic group and $X$ is the Cayley graph associated to the free product of $\Gamma$ and a free group (with respect to naturally defined generators), then $\isom(X)$ is a locally compact group which is not discrete (nor an extension of such), acting cocompactly on $X$, and Theorem~\ref{thm_superrigidity} applies. Thus this theorem covers completely new situations in terms of the target (and source) groups being involved. It is natural to expect that the cocompactness assumption of~(ii) is redundant (but then amenability of $H$ is to be replaced by elementarity).

\medskip

\noindent
\textbf{1.6.} The results of~\cite{Monod-Shalom1} together with the above Theorem~\ref{thm_coho} provide us with a large class of groups that are in~$\mathcal{C}_\mathrm{reg}$ because of various instances of negative curvature phenomena. It is therefore natural to study the class~$\mathcal{C}_\mathrm{reg}$ more thoroughly, and particularly to ask: \itshape Is there a characterisation of~$\mathcal{C}_\mathrm{reg}$ in geometric terms? Is being in~$\mathcal{C}_\mathrm{reg}$ a geometric property?\upshape

\smallskip

On the one hand, it would be interesting to find geometric restrictions on this class~--- in addition to the algebraic restrictions that we know~\cite{Monod-Shalom1}. On the other, one could try to include more groups in~$\mathcal{C}_\mathrm{reg}$ by searching amongst further relatives in the family of ``negative curvature'': for instance groups with a non-elementary proper isometric action on some general geodesic Gromov-hyperbolic metric space. Here we have in mind the possibility of non-proper spaces aswell, so that the notion of ``proper action'' should be that of \emph{metrically proper}, or maybe even weaker, as in~\cite{Bestvina-Fujiwara}. This would include mapping class groups, since they act on the curve complex, which is not locally finite but hyperbolic~\cite{Masur-Minsky}. A positive indication in that direction is the fact~\cite{Bestvina-Fujiwara} that groups with such actions have infinite dimensional $\hb^2(-,{\bf R})$; indeed, in concrete situations, it is often possible to recover (infinitely many) classes in $\hb^2(-,{\bf R})$ from the constructions that we use for $\hb^2(-,\ell^2(-))$. For more ``negative curvature'' aspects of the mapping class groups, compare also~\cite{Farb-Lubotzky-Minsky}.

\medskip

Finally, we ask: \emph{What (else) can be said about the class of all groups ME to some non-elementary hyperbolic group?}

In fact, as pointed out in~\cite{Monod-Shalom2}, the free group on two generators in already ME to the uncoutable class of all free products of any finite collection of (at least two) countable amenable groups.

%%%%%%%%%%%%%%%%%%%%%%%%%%%%%%%%%%%%%%%%%%%%%%%%%%%%%%%%%%%%%%%%%%%% 
\section{Notation}
\label{sec_notation}%

For any set $X$ we write $\ell^p(X)$ for the space of $p$\ti summable functions ($1\leq p< \infty$) and $\ell^\infty(X)$ for bounded functions; $C(X)$, $\coo(X)$
stand for the space of all, respectively all finitely supported functions; $\co(X)$ for those vanishing at infinity (\emph{i.e.} along the filtre of cofinite subsets). The Banach spaces $\ell^p(X)$ are endowed with the $p$\ti norms~$|\!|\cdot{|\!|}_p$, whilst the \weak topology on~$\ell^p(X)$ refers to the canonical isomorphisms $\ell^p(X)\cong \ell^{p/(p-1)}(X)^*$ for $1<p<\infty$, $\ell^\infty(X) \cong \ell^1(X)^*$ and $\ell^1(X)\cong \co(X)^*$.

We use the Serre's notation for graphs, so that a graph $\G=(\go,E)$ consists of a set of vertices $\go$, a set of edges $E$, an involution $e\mapsto\bar e$ of $E$ and initial/terminal maps $E\to\go$ (which we denote $e\mapsto e_-$, $e\mapsto e_+$) subject to the usual conditions. We say that $\G$ has bounded valency if there is a bound on the number of edges issuing from any vertex. As usual, we endow the geometric realization of $\G$ with the path metric $d$ of unit edge length and often abuse notation in writing $\G$ for the resulting space. Likewise, paths are simplicial paths aswell as path in the realization. The ``distance'' $d$ between subsets of $\G$ is the infimum of the distance between their points.

Elements of $\coo(\go)$ and $\coo(E)$ are also referred to as $0$\ti chains and $1$\ti chains. The boundary $\partial f$ of a $1$\ti chain $f$ is the $0$\ti chain $\partial f(v)= \sum_{e_+ =v} f(e) - \sum_{e_- =v} f(e)$; this definition extends to $C(E)$ if $\G$ is locally finite. Any (vector-valued) function on vertices or edges extends by linearity to a function on chains, and we will use the same notation (considering chains as formal sums). We also use the notation $\big< f,e\big>$ for the value of $f\in C(E)$ at the edge $e$, and $supp(f)$ is the closure of the union of the edges where $f$ does not vanish. Throughout the paper, $\G=(\go,E)$ is a hyperbolic graph of bounded valency; consequently, for every $R>0$ there is a uniform bound on the number of edges or vertices in any ball of radius $R$. For example, the Cayley graph of a hyperbolic group with respect to any finite set of generators would be such.

Let $(X,d)$ be a metric space. Denote by $B(x,r)$, respectively $\overline{B}(x,r)$, the open, respectively closed, {\sf $r$\ti ball} with centre~$x$, and by $N(S,r)$ the {\sf closed $r$\ti neighbourhood} of the subset~$S$ in~$X$. One calls $X$ {\sf proper} if all closed balls are compact. A background reference for {\sf hyperbolic} metric spaces is~\cite{Ghys-Harpe}. Unless otherwise stated, $H=\isom(X)$ is the group of isometries of $X$; in the graph case $X=\G$, it coincides with the group of automorphisms and all spaces of functions on $\G$ defined above are endowed with the natural $H$\ti representation. For any proper metric space $X$, the compact-open topology turns $H$ into a locally compact topological group and its action on $X$ is proper. One says that $X$ has {\sf at most exponential growth} if there is a constant $\alpha$ such that any ball of any radius $r$ contains at most $\alpha^r$ disjoint balls or radius one.

Hyperbolicity of~$X$ guarantees the existence of a constant $\de\ge 0$ such that all the geodesic triangles in $X$ are {\sf $\de$\ti fine} in the following sense: If $a,b,c$ are points of $X$ and $[a,b]$, $[b,c]$, and $[c,a]$ are geodesics from $a$ to $b$, from $b$ to $c$, and from $c$ to $a$, respectively, if moreover points $\bar{a}\in [b,c]$, $v,\bar{c}\in [a,b]$, $w,\bar{b}\in [a,c]$ satisfy
$$d(b,\bar{c})=d(b,\bar{a}),\quad d(c,\bar{a})=d(c,\bar{b}),\quad
d(a,v)=d(a,w)\le d(a,\bar{c})=d(a,\bar{b}),$$
then $d(v,w)\le \de$. We call $\{\bar{a},\bar{b},\bar{c}\}$ a triple {\sf inscribed} in the triangle $\{a,b,c\}$. For the rest of the paper, $\de$ denotes a positive integer, depending only on~$X$, such that all the geodesic triangles in~$X$ are $\de$\ti fine.

We write $\p X$ for the Gromov boundary and $\partial^n X\subseteq (\p X)^n$ for the subset of $n$\ti tuples of pairwise distinct points. The {\sf Gromov product} of $a,b\in X$ with respect to a basepoint $x_0\in X$ is
$$(a|b)_{x_0}:= \frac{1}{2}\big( d(a,x_0)+ d(b,x_0)- d(a,b)\big).$$
There is a natural extension of this product to $a,b\in\overline{X}:= X\sqcup \p X$; this extension is not necessarily continuous and in general one can make various choices of such extensions that differ by a bounded amount. The sets $\{x\in \overline{X} : (x|\xi)_{x_0}>r\}$ form a basis of neighbourhoods of $\xi\in\p X$ in~$\overline{X}$. The various extensions of $(\cdot|\cdot)_{x_0}$ induce the same topology on~$\p X$. In the graph case, we call a subgraph $U\se\G$ a {\sf graph neighbourhood} of $\xi\in\p\G$ if $U$ is the restriction to $\G$ of a neighbourhood of $\xi$ in $\GG$.

A subgroup of a hyperbolic group is said to be {\sf elementary} if it is virtually cyclic; elementary subgroups of hyperbolic groups are exactly the amenable subgroups (see \emph{e.g.} Thm.~38, p.~21 in \cite{Ghys-Harpe}). This notion has a natural generalisation in view of the following fact:

\begin{ppp}
\label{prop_def_elem}%
Let $X$ be a hyperbolic proper geodesic metric space and $H<\isom(X)$ any subgroup. The following are equivalent:

\begin{itemize}
\item[(i)] $H$ fixes a probability measure on $\p X$.

\item[(ii)] $X$ preserves a compact set in $X$ or a point in $\p X$ or a pair of points in $\p X$.
\end{itemize}

\noindent
Moreover, if $X$ has at most exponential growth (\emph{e.g.} if it is a graph of bounded valency or if $\isom(X)$ acts cocompactly), these are equivalent to

\begin{itemize}
\item[(iii)] The closure of $H$ in $\isom(X)$ is an amenable group.
\end{itemize}
\end{ppp}

For a proof see Section~\ref{sec_general_spaces}; we say that $H$ is {\sf elementary} if it satisfies the equivalent conditions~(i) and~(ii).

%%%%%%%%%%%%%%%%%%%%%%%%%%%%%%%%%%%%%%%%%%%%%%%%%%%%%%%%%%%%%%%%%%%%
\section{The Homological Bicombing and Exponential Convergence}

Let $\G$ be a hyperbolic graph of bounded valency.

In this section, we elaborate on the construction of a homological bicombing from~\cite{Mineyev01}, and establish its further properties which we need for its extension to infinity in order to obtain the ideal bicombing. For every $a,b\in\go$, let $S_{a,b}$ be the (finite) set of all geodesic paths from $a$ to $b$. We fix once and for all a choice $p[a,b]\in S_{a,b}$ for every vertices $a,b\in\go$. We will view $p[a,b]$ both as a path and as a $1$\ti chain. By $p[a,b](r)$ we mean the point on $p[a,b]$ (viewed as a path) at distance $r$ from~$a$. Contrary to the situation in~\cite{Mineyev01}, we cannot assume $p$ to be equivariant; we have however an $H$\ti equivariant map $p':\go\times \go\to \coo(E)$ defined by
$$p'[a,b] := |S_{a,b}|^{-1}\sum_{s\in S_{a,b}} s.$$
(So in particular $p'[a,b](r)$ will be viewed as a $0$\ti chain.) This $p'$ will be used instead of $p$ in the construction in order to get equivariance; however, the non-equivariant choice of geodesics $p$ fixed above will be used as auxilliary tool in some proofs.

\begin{ppp}[Proposition~7 in~\cite{Mineyev01}]
\label{bicombing}%
There exists a function $\bar{f}:\go\times\go\to \coo(\go)$ with the following properties.
\begin{itemize}
\item [(i)] $\bar{f}(b,a)$ is a convex combination, that is its coefficients are non-negative and sum up to~$1$.
\item [(ii)] If $d(a,b)\ge 10\de$, then $supp \bar{f}(b,a)\se \overline{B}(s(10\de), 8\de)$ for every $s\in S_{b,a}$.
\item [(iii)] If $d(a,b)\le 10\de$, then $\fff(b,a)=a$.
\item [(iv)] $\bar{f}$ is $H$\ti equivariant, \emph{i.e.} $\bar{f}(h b,h a)=h\,\bar{f}(b,a)$ for any $a,b\in\go$ and $h\in H$.
\item [(v)] There exist constants $L\in [0,\infty)$ and $\la\in[0,1)$ such that, for any $a,a',b\in\go$,
$$|\!| \bar{f}(b,a)-\bar{f}(b,a'){|\!|}_1\le L\la^{(a|a')_b}.$$
\item [(vi)] There exists a constant $\la'\in[0,1)$ such that if $a,b,b'\in\go$ satisfy $(a|b')_b\le 10\de$ and $(a|b)_{b'}\le 10\de$, then
$$|\!|\bar{f}(b,a)-\bar{f}(b',a){|\!|}_1\le 2\la'.$$
\item [(vii)] Let $a,b,c\in\go$, $s\in S_{a,b}$, and let $c\in N(s, 9\de)$. Then $supp\, \bar{f}(c,a)\se N(s, 9\de)$. 
\end{itemize}
\end{ppp}

\begin{proof}[About the proof]
\textbf{(a)}~The construction of $\fff(a,b)$ in~\cite{Mineyev01} made use of an equivariant choice of geodesic paths $p[a,b]$. If we replace this $p$ with the above $p'$, the construction and arguments in~\cite{Mineyev01} can easily be adapted to yield the above $\bar{f}$; it is important here that there is a uniform bound on the number of vertices in any ball of a given radius.

\smallskip

\textbf{(b)}~In addition, we slightly redefine the function $\fff$ given in~\cite[Proposition~7]{Mineyev01} by imposing the above condition~(iii).
\end{proof}

\subsection{The homological \boldmath bicombings~$q'$ and~$q$ \unboldmath}

Taking into account that we use $p'$ instead of a choice of geodesics, the homological bicombing $q':\go\times\go\to \coo(E)$
constructed in~\cite{Mineyev01} can now be defined as follows.
Set $q'[a,a]:=0$; then, inductively on $d(a,b)$,
$$q'[a,b]:= q'[a,\bar{f}(b,a)]+ p'[\bar{f}(b,a),b].$$
We recall that the word (homological) {\sf bicombing} refers to the property
\begin{equation}
\label{eq_bic}%
\p q'[a,b] = b-a\kern1cm \forall\,a,b\in\go.
\end{equation}

\begin{ppp}[Proposition~8 in~\cite{Mineyev01}]
\label{q'}%
The homological bicombing $q'$ above satisfies the following conditions.
\begin{itemize}
\item [(i)] $q'$ is $\Gamma$\ti equivariant, \emph{i.e.} $\gamma\cdot q'[a,b]=q'[\gamma\cdot a, \gamma\cdot b]$
for $a,b\in\go$ and $\gamma\in \gamma$.
\item [(ii)] $q'$ is quasigeodesic. More precisely, for all $a,b\in\go$, 
the support of $q'[a,b]$ lies in the $27\de$\ti neighbourhood of any geodesic from $a$ to $b$; moreover $|\!|q'[a,b] {|\!|}_1\le 18\de d(a,b)$.
\item [(iii)] There exist constants $M\in[0,\infty)$ and $N\in[0, \infty)$ 
such that,
for all $a,b,c\in\go$,
$$|\!| q'[a,b]- q'[a,c] {|\!|}_1  \le Md(b,c)+ N.$$\hfill\qedsymbol
\end{itemize}
\end{ppp}

The homological bicombing $q$ is now defined as follows:
$$q[a,b]:= \frac{1}{2} \big(q'[a,b]- q'[b,a] \big).$$
It was shown in~\cite{Mineyev01} that $q$ has bounded area, that is,
\begin{equation}
\label{eq_b_area}%
\sup_{a,b,c\in \go} |\!| q[a,b]+ q[b,c]+ q[c,a] {|\!|}_1 < \infty.
\end{equation}

\subsection{More properties of $q'$}

\begin{ppp}
\label{2003de}%
For all vertices $a,b$ and any edge $e$ in $\G$,\ $\big|\big<q'[a,b], 
e\big>\big|\le 2003\,\de^2$.
\end{ppp}

\begin{proof} 
Fix a vertex~$a$ and an edge $e$ in~$\G$.

{\bf (a)} First assume $d(a,e_-)\ge 50\de$.

The definition of $q'$ and Proposition~\ref{bicombing}(vii) 
imply that the support of $q'[a,b]$ lies in the $27\de$\ti neighbourhood
of any geodesic between $a$ and $b$. So if a vertex $b$ satisfies 
$d(a,b)\le d(a,e_-)-28\de$,
then the support of $q'[a,b]$ is disjoint from $e$ and
$\left|\left<q'[a,b], e\right>\right|=0$.

Now assume
\begin{equation}
\label{28de}%
d(a,e_-)-28\de \le d(a,b)\le d(a,e_-)+20\de.
\end{equation}
In this case we prove the inequality
\begin{equation}
\label{18de}%
\big|\big<q'[a,b], e\big>\big|\le 18\de\big(d(a,b)-d(a,e_-)+46\de\big)
\end{equation}
inductively on $d(a,b)$.

If $b$ is such that $d(a,e_-)-46\de\le d(a,b)\le d(a,e_-)-28\de$,
then by the argument above,
$$\big< q'[a,b], e\big> =0\le 18\de\big(d(a,b)-d(a,e_-)+46\de\big),$$
which provides for the basis of induction.
For the inductive step, if $b$ satisfies~(\ref{28de}), then 
each vertex~$x$ in the support of $\bar{f}(b,a)$ satisfies
$$d(a,e_-)-46\de\le d(a,x)\le d(a,b)-1,$$
so the induction hypothesis applies to each such $x$ giving
\begin{eqnarray*}
&&\left|\left< q'[a,b],e\right>\right|= \left|\left<q'[a,\bar{f}(b,a)]+ 
p'[\bar{f}(b,a),b],e\right>\right|
\le\left|\left<q'[a,\bar{f}(b,a)],e\right>\right|+ 
\left|\left<p'[\bar{f}(b,a),b],e\right>\right|\\
&& \le 18\de\big( d(a,b) -1 -d(a,e_-)+ 46\de\big)+ 
 |\!| p'[\bar{f}(b,a),b]{|\!|}_1\\
&&\le 18\de\big( d(a,b) -1 -d(a,e_-)+ 46\de\big)+ 18\de= 18\de\big( 
d(a,b)-d(a,e_-)+ 46\de\big).
\end{eqnarray*}

This proves~(\ref{18de}) in the case when~(\ref{28de}) holds, so this gives 
the bound
\begin{equation}
\label{2003}%
\big|\big<q'[a,b],e\big>\big|\le
18\de\big(d(a,b)-d(a,e_-)+46\de\big)\le 18\de\big( 20\de+46\de\big)\le 
2003\,\de^2.
\end{equation}
Finally, assume $d(a,b)\ge d(a,e_-)+20\de$. Inductively on $d(a,b) $ we 
show that
the above bound is preserved. The basis of induction is provided 
by~(\ref{2003}), and 
for the inductive step note that since the support of $p'[\bar{f}(b,a),b]$ 
lies in the ball $\overline{B}(b,18\de)$, it is disjoint from $e$, so
\begin{eqnarray*}
&&\big|\big<q'[a,b],e\big>\big|= \big|\big<q'[a,\bar{f}(b,a)]+ 
p'[\bar{f}(b,a),b],e\big>\big|=
\big|\big<q'[a,\bar{f}(b,a)], e\big>\big|\le 2003\,\de^2.
\end{eqnarray*}

{\bf (b)}
Now assume $d(a,e_-)\le 50\de$.
If $d(a,b)\le d(a,e_-)+20\de$, then $d(a,b)\le 50\de+ 20\de=70\de$ and
Proposition~\ref{q'}(ii) gives the same bound as in~(\ref{2003}):
$$\left|\left<q'[a,b],e\right>\right|\le |\!|q'[a,b]{|\!|}_1\le 
18\de\cdot 48\de\le 2003\,\de^2.$$
In the case $d(a,b)\ge d(a,e_-)+20\de$ this bound is preserved, similarly to (a) 
above.
This finishes the proof of Proposition~\ref{2003de}.
\end{proof}

\begin{lem}
\label{PQ}
There exist constants $Q\in [0,\infty)$ and $\sigma_0\in [0,1)$ such that 
for all
$a,b,b'\in\Gamma$ and any edge $e$ in $\G$, if $d(b,b')\le 56\de$ then
$$\left| \left< q'[a,b]-q'[a,b'],e\right> \right|\le
(d(b,b')+ Q)\, \sigma_0^{d(e_-,b)+d(e_-,b')}.$$
\end{lem}

\begin{proof} 
First we define the constants.
Let $\la'\in[0,1)$ be the constant from Proposition~\ref{bicombing}(vi).
Choose $\sigma_0\in[0,1)$ close enough to~1 so that
\begin{equation}
\label{sigma_0-1}%
1-\la' \sigma_0^{-148\de}> 0.
\end{equation}

Next, choose $Q$ large enough so that
\begin{equation}
\label{Q1}%
2\cdot 2003\,\de^2\le Q \sigma_0^{312\de} \mbox{\qquad and\qquad}
\la'\cdot 56\de\sigma_0^{-148\de}\le Q(1-\la' \sigma_0^{-148\de}).
\end{equation}

Then make $\sigma_0$ closer to~1, if needed, so that
\begin{equation}
\label{sigma_0-2}%
(56\de+ Q)(\sigma_0^{-74\de}-1)\le \sigma_0^{-74\de}.
\end{equation}
Note that (\ref{sigma_0-1}) and (\ref{Q1}) still hold for this new $\sigma_0$.
The second inequality in~(\ref{Q1}) rewrites as
\begin{equation}
\label{Q2}%
\la' (56\de+ Q)\, \sigma_0^{-148\de}\le Q.
\end{equation}
Each of these constants depends only on~$\G$.
Now we proceed with the proof.

Fix a vertex~$a$ and an edge~$e$ in~$\G$.
We will use induction on $d(a,b)+d(a,b')$ for vertices $b$ and~$b'$.

If $e$ is not in the $28\de$\ti neighbourhood of either $p[a,b]$ or $p[a,b']$, then
by Proposition~\ref{q'}(ii),
$\left|\left< q'[a,b]-q'[a,b'],e\right> \right|=0$ and the result follows.
If~$e$ is in the $28\de$-neighbourhood of~$p[a,b]$ but not of $p[a,b']$, then
$e_-$ is $\de$-close to $p[b,b']$ and by~(\ref{Q1}),
\begin{eqnarray*}
&&|\left<q'[a,b]-q'[a,b'], e\right>|\le |\left<q'[a,b],e\right> |\le\\
&&\ \ \le 2003\de^2 \le Q\sigma_0^{312\de}\le Q\sigma_0^{d(b,b')+2\de}\le Q\sigma_0^{d(e_-,b)+d(e_-,b')}.
\end{eqnarray*}
Similarly for~$b$ and~$b'$ interchanged, so from now on we can assume that $e$
is in the $28\de$\ti neighbourhoods of $p[a,b]$ and $p[a,b']$.

First assume that
\begin{equation}
\label{2dae}%
d(a,b)+d(a,b')\le 2 d(a,e_-)+100\de.
\end{equation}
If $d(a,b)> d(a,e_-)+100\de$ were true, then
$d(a,b')> d(a,b)-100\de\ge d(a,e_-)$ and 
$d(a,b)+ d(a,b')> 2\,d(a,e_-)+ 100\de$,
which is a contradiction. Therefore, $d(a,b)\le d(a,e_-)+100\de$. Similarly,
we obtain $d(a,b')\le d(a,e_-)+100\de$.

Let $u$ be a vertex on the geodesic $p[a,b']$ closest to $e_-$. Then
\begin{eqnarray*}
&&d(e_-,b')\le d(e_-,u)+ d(u,b')\le 28\de+ d(a,b')- d(a,u)\\ 
&& \le 28\de+ d(a,b')- (d(a,e_-)- 28\de)\le d(a,b')- d(a,e_-)+ 56\de\le 100\de+56\de= 156\de.
\end{eqnarray*}
Similarly, $d(e_-,b)\le 156\de$.
Then using Proposition~\ref{2003de} and (\ref{Q1}),
\begin{eqnarray*}
&&\left| \left< q'[a,b]-q'[a,b'],e\right> \right|\le
 \left| \left< q'[a,b],e\right> \right|+ \left| \left< q'[a,b'],e\right> \right|
 \le 2\cdot 2003\,\de^2\\
&& \le Q\, \sigma_0^{312\de}
 \le Q\, \sigma_0^{d(e_-,b)+d(e_-,b')}\le (d(b,b')+ Q)\, \sigma_0^{d(e_-,b)+d(e_-,b')},
\end{eqnarray*}
which proves the desired inequality assuming~(\ref{2dae}).

From now on we can assume $d(a,b)+d(a,b')\ge 2 d(a,e_-)+100\de$.
Since $d(b,b')\le 56\de$, this implies
\begin{equation}
\label{dab1}%
   d(a,b)\ge d(a,e_-)+ 20\de \mbox{\qquad and\qquad} d(a,b')\ge d(a,e_-)+ 20\de
\end{equation}
(by contradiction). Consider the following two cases.

{\bf Case~1.} $(a|b')_{b}> 10\de$ or $(a|b)_{b'}> 10\de$ (see Fig.~\ref{case1}).\\
Without loss of generality, $(a|b')_{b}> 10\de$, the other
case being similar. Let $x$ be an arbitrary vertex in the support of $\bar{f}(b',a)$
and $v:=p[b,a](10\de)$.
Then by Proposition~\ref{bicombing}(ii) we have
\begin{eqnarray}
\label{ex}%
&&d(e_-,x)\ge d(u,v)-8\de-28\de\ge d(u,b)-10\de-8\de-28\de\\
\nonumber&&= d(u,b)-46\de\ge d(e_-,b)-28\de-46\de= d(e_-,b)-74\de.
\end{eqnarray}

\setlength{\unitlength}{1 cm} 
\begin{figure}[h]
  \begin{center}
   \begin{picture}(8,4)
%\input{fig/sh_case1}
%\graphpaper[1](0,0)(8,4)

\put(0,2){\circle*{0.1}}
\put(-.4,1.6){\footnotesize$a$}

\put(7,4){\circle*{0.1}}
\put(7.2,4){\footnotesize$b$}

\put(6,0){\circle*{0.1}}
\put(6.2,0){\footnotesize$b'$}

%from a to b
\qbezier(0,2)(4.4,2)(5.1,2.3)
\qbezier(5.1,2.3)(5.6,2.5)(7,4)
\put(1.7,2.3){\footnotesize$p[b,a]$}
%from a to b'
\qbezier(0,2)(4.1,2)(4.8,1.63)
\qbezier(4.8,1.63)(5.2,1.45)(6,0)
\put(1.6,1.6){\footnotesize$p[b',a]$}

%from b to b'
\qbezier(7,4)(5,2)(6,0)
\put(6.7,3.4){\tiny$10\de$}
%inscribed triangle
\put(5.07,2.28){\circle*{0.1}}
\put(5.7,1.7){\circle*{0.1}}
\put(5,1.49){\circle*{0.1}}
%v
\put(6.04,3.02){\circle*{0.1}}
\put(5.7,2.83){\footnotesize$v$}
%v'
\put(6.16,2.95){\circle*{0.1}}

\put(5.1,3.04){\circle*{0.1}}
\put(4.65,3.04){\footnotesize$x$}

%from v to x
\qbezier(6.04,3.02)(5.57,3.03)(5.1,3.04)
\put(5.4,3.2){\tiny$\le\! 8\de$}
%from v to v'
\qbezier(6.04,3.02)(6.1,2.98)(6.16,2.95)

\put(4,1.82){\circle*{0.1}}
\put(3.62,1.6){\footnotesize$u$}

\put(3.7,0.8){\circle*{0.1}}
\put(3.4,0.6){\footnotesize$e_-$}

%from u to e_-
\qbezier(4,1.82)(3.85,1.31)(3.7,0.8)
\put(3.9,1.2){\tiny$\le\! 28\de$}

   \end{picture}
  \end{center}
 \caption{\label{case1} Case~1.}
\end{figure}

Inequality~(\ref{sigma_0-2}) implies
\begin{equation*}
\big( d(b,b')+ Q\big) (\sigma_0^{-74\de}-1)\le 
\big( 56\de+ Q\big) (\sigma_0^{-74\de}-1)\le \sigma_0^{-74\de},
\end{equation*}

which is equivalent to
\begin{equation}
\label{Pbb'}%
\big( d(b,b')- 1+ Q\big) \sigma_0^{-74\de}\le d(b,b')+ Q.
\end{equation}

It is easy to check that
\begin{equation}
\label{bx}%
d(b',x)\le d(b,b')-1\le 56\de\mbox{\qquad and\qquad} d(a,x)+d(a,b')< d(a,b)+d(a,b'),
\end{equation}
therefore
the induction hypotheses are satisfied for the vertices~$x$ and~$b'$. Using the induction
hypotheses, (\ref{ex}), (\ref{bx}) and~(\ref{Pbb'}), we have
\begin{eqnarray}
\label{q'abq'ax}%
&&\left|\left< q'[a,b']-q'[a,x],e\right>\right|\le \big( d(b',x)+Q\big)\,\sigma_0^{d(e_-,x)+d(e_-,b')}\\
\nonumber &&\le \big( d(b,b')- 1+ Q\big)\, \sigma_0^{d(e_-,b)- 74\de+ d(e_-,b')}
  = \big( d(b,b')- 1+ Q\big)\, \sigma_0^{-74\de} \sigma_0^{d(e_-,b)+d(e_-,b')}\\
\nonumber &&\le \big( d(b,b')+ Q\big)\, \sigma_0^{d(e_-,b)+d(e_-,b')}.
\end{eqnarray}
By (\ref{dab1}), $e$ is disjoint from $supp(p'[x,b])$ for each
$x\in supp(\fff(b,a))$, therefore~$e$ is disjoint from $supp(p'[\fff(b,a),b])$.
Since $q'[a,x]$ is linear in $x$, (\ref{q'abq'ax}) implies the desired inequality:
\begin{eqnarray*}
&& \left|\left< q'[a,b]-q'[a,b'],e\right>\right|
 = \left|\left< q'[a,\fff(b,a)]+ p'[\fff(b,a),b]- q'[a,b'],e\right>\right|\\
&& = \left|\left< q'[a,\fff(b,a)]- q'[a,b'],e\right>\right|\le 
  \big( d(b,b')+ Q\big)\, \sigma_0^{d(e_-,b)+d(e_-,b')}.
\end{eqnarray*}

{\bf Case~2.} $(a|b')_b\le 10\de$ and $(a|b)_{b'}\le 10\de$.

In particular, $d(b,b')= (a|b')_b+ (a|b)_{b'}\le 20\de$.
Since
$$d(a,b)+d(a,b')\ge 2d(a,e_-)+100\de> 40\de,$$
one can easily see that $d(a,b)>10\de$ and $d(a,b')>10\de$.

The 0-chain $\fff(b,a)-\fff(b',a)$ has the form $f_+-f_-$,
where $f_+$ and $f_-$ are 0-chains with
non-negative coefficients and disjoint supports.
By Proposition~\ref{bicombing}(6),
$$|f_+|_1+|f_-|_1= |f_+ - f_-|_1=\big|\fff(b,a)-\fff(b',a)\big|_1\le 2\la',$$
where $\la'\in[0,1)$ is independent of $a,b,b'$.
The coefficients of the 0-chain $f_+ - f_-=\fff(b,a)-\fff(b',a)$
sum up to 0, because $\fff(b,a)$ and $\fff(b',a)$ are convex
combinations. It follows that
\begin{equation}
\label{f+-}%
|f_+|_1=|f_-|_1\le\la'.
\end{equation}
Also,
$$supp\,f_+\se
  supp\,\fff(b,a)\se
  \overline{B}(p[b,a](10\de),8\de)$$
and
$$supp\,f_-\se
  supp\,\fff(b',a)\se
  \overline{B}(p[b',a](10\de),8\de),$$
hence, by the hypotheses of Case~2, for all $x\in supp\,f_+$ and $x'\in supp\,f_-$,
$$d(x,x')\le d(x,b)+ d(b,b')+ d(b',x')\le
  18\de+ 20\de+ 18\de= 56\de.$$
Also $d(a,x)+d(a,x')< d(a,b)+d(a,b')$, so all such pairs~$x$ and~$x'$ satisfy the induction hypotheses.
Similarly to~(\ref{ex}), we obtain
$$d(e_-,x)\ge d(e_-,b)-74\de \mbox{\quad and\quad} d(e_-,x')\ge d(e_-,b')-74\de,$$
then using the induction hypotheses,
\begin{eqnarray}
\label{q'axq'ax'}%
&& \left|\left<q'[a,x]-q'[a,x'],e\right>\right|\le
  \big( d(x,x')+ Q\big)\, \sigma_0^{d(e_-,x)+ d(e_-,x')}\\
\nonumber && \le (56\de+Q)\, \sigma_0^{d(e_-,b)-74\de+ d(e_-,b')- 74\de}=
  (56\de+Q)\, \sigma_0^{-148\de} \sigma_0^{d(e_-,b)+ d(e_-,b')}.
\end{eqnarray}

Using~(\ref{dab1}) we see that $e$ is disjoint from $p'[\fff(b,a),b]$ and $p'[\fff(b',a),b']$, then
by linearity of $q'[a,x]$ in $x$, (\ref{f+-}), (\ref{q'axq'ax'}), and (\ref{Q2}),
\begin{eqnarray*}
&& \left|\left<q'[a,b]-q'[a,b'],e\right>\right|\\
&& \le \left|\left<q'[a,\fff(b,a)]+ p'[\fff(b,a),b]- q'[a,\fff(b',a)]- p'[\fff(b',a),b'],e\right>\right|\\
&&\le \left|\left<q'[a,\fff(b,a)]-q'[a,\fff(b',a)],e\right>\right|\le
  \left|\left<q'[a,f_+]-q'[a,f_-],e\right>\right|\\
&& \le \la'\cdot (56\de+Q) \sigma_0^{-148\de} \sigma_0^{d(e_-,b)+ d(e_-,b')}
   \le Q\,\sigma_0^{d(e_-,b)+ d(e_-,b')}\\
&& \le \big( d(b,b')+Q\big) \sigma_0^{d(e_-,b)+ d(e_-,b')},
\end{eqnarray*}
which is the desired inequality. This finishes the proof of Lemma~\ref{PQ}.
\end{proof}

\begin{ppp}
\label{p_bb'}%
There exist constants $S_1\in [0,\infty)$ and $\sigma_1\in [0,1)$
such that, for all vertices $a,b,b'$ with $d(b,b')\le 1$ and any edge $e$ in $\G$,
$$\left| \left<q'[a,b]- q'[a,b'],e \right>\right|\le S_1\,\sigma_1^{d(e_-,b)}.$$
\end{ppp}

\begin{proof}
Immediate from Lemma~\ref{PQ}: since $d(b,b')\le 56\de$, we have
\begin{eqnarray*}
&&\left| \left< q'[a,b]-q'[a,b'],e\right> \right|\le
(d(b,b')+ Q)\, \sigma_0^{d(e_-,b)+d(e_-,b')}\\
&&\le (1+Q)\, \sigma_0^{d(e_-,b)+d(e_-,b)-1}=
  (1+Q)\sigma_0^{-1} (\sigma_0^2)^{d(e_-,b)},
\end{eqnarray*}
so we set $S_1:=(1+Q)\sigma_0^{-1}$ and $\sigma_1:=\sigma_0^2$.
\end{proof}

Now we would like to have a similar result with the roles of $a$'s and $b$'s
interchanged (Proposition~\ref{p_aa'}). Since the construction of $q'[a,b]$ is not symmetric
in~$a$ and~$b$, a different argument is required. We start with the following lemma.

\begin{lem}
\label{Kla}%
There exist constants $K\in[0,\infty)$ and $\la\in [0,1)$ such that for all vertices
$a,a',b$ with $d(a,a')\le 1$ and any edge~$e$ in $\G$, the following conditions hold.
\begin{itemize}
\item [(i)]  If $d(a,b)\le d(a,e_-)- 30\de$, then $\left<q'[a,b]-q'[a',b],e\right>= 0$.
\item [(ii)] If $d(a,b)\ge d(a,e_-)- 60\de$, then
  $$\left|\left<q'[a,b]-q'[a',b],e\right>\right|\le
  K\left(\la^{d(a,e_-)- 60\de}+ \la^{d(a,e_-)- 60\de+1}+ ...+ \la^{d(a,b)}\right).$$
\end{itemize}
\end{lem}

\begin{proof}
Let $\la\in[0,1)$ and $L\in[0,\infty)$ be the constants
from Proposition~\ref{bicombing}(v), and choose $K$ large enough so that
\begin{equation}
\label{K}%
2\cdot 2003\,\de^2\le K\mbox{\qquad and\qquad} 2L\la^{-1}\cdot 2003\,\de^2\le K.
\end{equation}

Fix vertices $a$ and $a'$ with $d(a,a')\le 1$ and an edge~$e$ in~$\G$.

{\bf (i)} Assume that $b$ satisfies $d(a,b)\le d(a,e_-)-30\de$. By Proposition~\ref{q'}(ii),
$q'[a,b]$ and $q'[a',b]$ lie in the $27\de$\ti neighbourhoods of the geodesics
$p[a,b]$ and $p[a,b']$, respectively. It is easy to see that~$e$ is
disjoint from these neighbourhoods, so (i) follows.

{\bf (ii)} Assume that $b$ satisfies $d(a,b)\ge d(a,e_-)-60\de$.

{\bf (a)} First suppose that $d(a,e_-)\le 60\de$. Then by Proposition~\ref{2003de}
and~(\ref{K}),
\begin{eqnarray*}
&& \left|\left< q'[a,b]- q'[a',b],e\right>\right|\le 
   \left|\left< q'[a,b],e\right>\right| + \left|\left< q'[a',b],e\right>\right|
  \le 2\cdot 2003\,\de^2\\
&& \le K \le K \la^{d(a,e_-)-60\de}
   \le K\big( \la^{d(a,e_-)-60\de}+ \la^{d(a,e_-)-60\de+1}+ ... +\la^{d(a,b)} \big).
\end{eqnarray*}

{\bf (b)} Now suppose that $d(a,e_-)\ge 60\de$.
We use induction on $d(a,b)$ to show (ii). Part~(i) provides the base of induction:
if $b$ satisfies
$$d(a,e_-)- 60\de\le d(a,b)\le d(a,e_-)-30\de,$$
then obviously
$$\left|\left< q'[a,b]-q'[a',b],e \right>\right|= 0\le
K\left(\la^{d(a,e_-)- 60\de}+ \la^{d(a,e_-)- 60\de+1}+ ...+ \la^{d(a,b)}\right).$$
For the inductive step, assume that $b$ satisfies $d(a,b)\ge d(a,e_-)-30\de$.
We have 
\begin{eqnarray}
\label{aa'b}%
&& (a|a')_b= \frac{1}{2} \big( d(a,b)+ d(a',b)- d(a,a') \big)\\
\nonumber && \ge \frac{1}{2} \big( d(a,b)+ d(a,b)-1-1 \big)= d(a,b)-1.
\end{eqnarray}
Since, for any 0-chain~$f$, $p'[f,b]$ is a convex combination of geodesics, and
any edge occurs in any geodesic at most once, it follows that
\begin{equation}
\label{eq_norm}%
\left|\left< p'[f,b],e\right>\right|\le |\!|f {|\!|}_1.
\end{equation}
Represent the 0-chain $\fff(b,a)-\fff(b,a')$ as $f_+ - f_-$ where
$f_+$ and $f_-$ have non-negative coefficients and disjoint supports.
Then there exists a 0-chain~$f_0$ with non-negative coefficients such that
$$\fff(b,a)=f_0+f_+,\qquad \fff(b,a')=f_0+f_-\mbox{,\qquad and\qquad}
  |\!|f_0{|\!|}_1\le 1.$$
Moreover, by Proposition~\ref{bicombing}(v),
\begin{equation}
\label{f-+}%
|\!|f_-{|\!|}_1= |\!|f_+{|\!|}_1=
   \frac{1}{2}|\!| \bar{f}(b,a)-\bar{f}(b,a'){|\!|}_1\le 
   \frac{1}{2} L\la^{(a|a')_b}.
\end{equation}
By Proposition~\ref{bicombing}(ii),
$$supp f_0= supp(\fff(b,a)) \cap supp(\fff(b,a'))\se \overline{B}(p[b,a](10\de),8\de)$$
therefore each vertex $x\in supp f_0$
satisfies $d(a,x)< d(a,b)$, so by the induction hypotheses,
\begin{eqnarray}
\label{axa'x}%
&& \left|\left< q'[a,x]-q'[a',x],e \right>\right|\le
     K\left(\la^{d(a,e_-)-60\de}+ ...+ \la^{d(a,x)}\right)\\
\nonumber && \le K\left(\la^{d(a,e_-)-60\de}+ ...+ \la^{d(a,b)-1}\right).
\end{eqnarray}
Using~(\ref{K}), (\ref{aa'b}), (\ref{axa'x}), (\ref{eq_norm}), (\ref{f-+}) we obtain
\begin{eqnarray*}
&& \left|\left< q'[a,b]-q'[a',b],e \right>\right| \\
&& = \left|\left< q'[a,\fff(b,a)]+ p'[\fff(b,a),b]- q'[a',\fff(b,a')]- p'[\fff(b,a'),b],e\right>\right|\\
&& = \left|\left< q'[a,\fff(b,a)]- q'[a',\fff(b,a')]+ p'[\fff(b,a)- \fff(b,a'),b],e\right>\right|\\
&& = \left|\left< q'[a,f_0+ f_+],e\right> - \left< q'[a,f_0+ f_-],e\right> +
    \left<p'[f_+- f_-,b],e\right>\right|\\
&& \le \left|\left< q'[a,f_0]- q'[a',f_0],e\right>\right|+
   \left|\left< q'[a,f_+],e\right>\right|+ \left|\left< q'[a,f_-],e\right>\right|+
   \left|\left< p'[f_+,b],e\right>\right|+ \left|\left< p'[f_-,b],e\right>\right|\\
&& \le |\!|f_0 {|\!|}_1\cdot K\left(\la^{d(a,e_-)-60\de}+ ...+ \la^{d(a,b)-1}\right)\\
&& \ \ \ + |\!|f_+ {|\!|}_1\cdot 2003\,\de^2+ |\!|f_- {|\!|}_1\cdot2003\,\de^2+
  |\!|f_+ {|\!|}_1+ |\!|f_- {|\!|}_1\\
&& \le K\left(\la^{d(a,e_-)-60\de}+ ...+ \la^{d(a,b)-1}\right)
 + |\!|f_+ {|\!|}_1\cdot 2003\,\de^2\cdot 4\\
&& \le K\left(\la^{d(a,e_-)-60\de}+ ...+ \la^{d(a,b)-1}\right)+
  \frac{1}{2} L\la^{(a|a')_b}\cdot 2003\,\de^2\cdot 4\\
&& \le K\left(\la^{d(a,e_-)-60\de}+ ...+ \la^{d(a,b)-1}\right)+
   2L \la^{d(a,b)-1}\cdot 2003\,\de^2\\
&& \le  K\left(\la^{d(a,e_-)-60\de}+ ...+ \la^{d(a,b)-1}\right)+
   K\la^{d(a,b)}\\
&& \le K\left(\la^{d(a,e_-)-60\de}+ ...+ \la^{d(a,b)}\right).
\end{eqnarray*}
This proves Lemma~\ref{Kla}.
\end{proof}

\begin{ppp}
\label{p_aa'}%
There exist constants $S_2\in[0,\infty)$ and $\sigma_2\in[0,1)$ such that,
for all vertices $a,a',b$ with $d(a,a')\le 1$ and any edge $e$ in $\G$,
$$\left|\left< q'[a,b]-q'[a',b], e\right>\right|\le S_2 \sigma_2^{d(a,e_-)}.$$
\end{ppp}

\begin{proof}
Immediate from Lemma~\ref{Kla}: if $d(a,b)\le d(a,e_-)-30\de$, then the statement
obviously follows, and if $d(a,b)\ge d(a,e_-)-60\de$, then
\begin{eqnarray*}
&& \left|\left< q'[a,b]- q'[a',b]\right>\right|\le
   K \left(\la^{d(a,e_-)-60\de}+ ...+ \la^{d(a,b)}\right)\\
&& \le K\la^{d(a,e_-)-60\de} \sum_{i=0}^\infty \la^i=
   \frac{K\la^{-60\de}}{1-\la} \la^{d(a,e_-)},
\end{eqnarray*}
so we set $S_2:= \frac{K\la^{-60\de}}{1-\la}$ and $\sigma_2:=\la$.
\end{proof}

Our goal is to combine Proposition~\ref{p_bb'} and Proposition~\ref{p_aa'}
into one stronger statement (Theorem~\ref{aa'bb'}).

Let $v,b,b'\in \go$, and let $\beta:J\to\G$
be a geodesic (\emph{i.e.} an isometric embedding of an interval~$J$) connecting~$b$ to~$b'$. 
Abusing the notation we will identify~$\be$ with the image of~$\be$.
Pick a distance-minimizing vertex $b_0\in\beta$, \emph{i.e.} such that $d(v,b_0)=d(v,\beta)$.
Let $\gamma=[v,b_0]$ be any geodesic between $a_0$ and $b_0$, then $\gamma$ is a shortest
geodesic connecting $v$ to~$\be$.
Choose the interval $J$ in $\rr$ containing 0 so that
$\beta(0)=b_0$, $\beta(d(b_0,b))=b$, $\beta(-d(b_0,b'))=b'$.

\begin{lem}
\label{ab6de}%
With the above notations, for all $j\in J$,
$$d(v,\beta(j))\ge d(v,\beta)+ |j|- 2\de.$$
\end{lem}

\begin{proof}
From symmetry, it suffices to show the lemma when $j\ge 0$.
\setlength{\unitlength}{1cm} 
\begin{figure}[h]
  \begin{center}
%\graphpaper[1](0,0)(9,5)
   \begin{picture}(9,4.5)

\put(0,0.5){\circle*{0.1}}
\put(0,0.7){\footnotesize$b'$}

\put(9,0){\circle*{0.1}}
\put(8.8,0.2){\footnotesize$b$}

% b to b'
\qbezier(0,0.5)(4.5,0.7)(9,0)

\put(4,4.6){\circle*{0.1}}
\put(4.2,4.4){\footnotesize$v$}

\put(3.8,0.51){\circle*{0.1}}
\put(2.5,0.1){\footnotesize$b_0=\be(0)$}

% v to b_0
\qbezier(4,4.6)(3.85,3.3)(3.8,0.51)

\put(7.5,0.2){\circle*{0.1}}
\put(7.1,0.5){\footnotesize$\be(j)$}

% v to b_0
\qbezier(4,4.6)(3.5,0.7)(7.5,0.2)
\put(3.5,2.3){\tiny$\gamma$}

%inscribe:
\put(4.79,1.36){\circle*{0.1}}
\put(5,1.5){\tiny$w_3$}

\put(4.4,0.48){\circle*{0.1}}
\put(4.3,0.1){\tiny$w_2$}

\put(3.82,1.1){\circle*{0.1}}
\put(3.3,1.2){\tiny$w_1$}
   \end{picture}
  \end{center}
\caption{Illustration for Lemma~\ref{ab6de}.}
\label{vbb'}%
\end{figure}
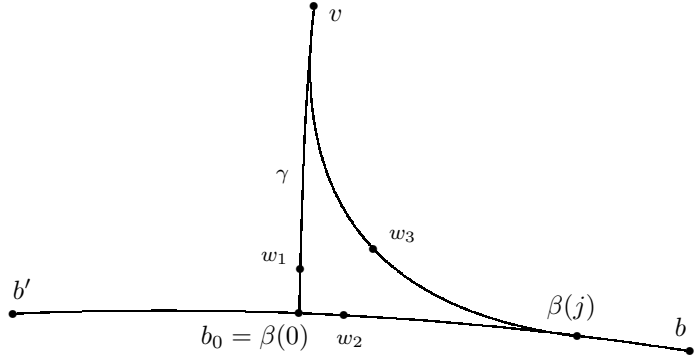
Inscribe a triple of points in the triangle
$\{v,b_0,\be(j)\}$ as shown on Fig.~\ref{vbb'}.
It is an easy exercise to see from the figure that
$d(b_0,w_2)= d(b_0,w_1)\le \de$, therefore
\begin{eqnarray*}
&& d(v,\be(j))= d(v,w_3)+ d(w_3,\be(j))= d(v,w_1)+ d(w_2,\be(j))\\
&& \ge \big( d(v,b_0)- \de\big)+ \big( d(b_0,\be(j)- \de \big)=
   d(v,\be)+ |j|- 2\de.
\end{eqnarray*}
\end{proof}

We recall that $q$ was defined by
\begin{equation}
\label{defi_q}%
q[a,b]= \frac{1}{2} \big(q'[a,b]- q'[b,a] \big)
\end{equation}
and proceed to state the exponential convergence:
\begin{ttt}
\label{aa'bb'}%
Let~$\Gamma$ be a hyperbolic group, $\G$ a graph with a proper cocompact $\Gamma$\ti action and $q$ the above bicombing. Then there exist constants $S\in[0,\infty)$ and $\sigma\in[0,1)$
such that, for all vertices $a,a',b,b'$ and any edge~$e$ in $\G$,
$$\left|\left<q[a,b]-q[a',b'],e\right>\right|\le S \left( \sigma^{(a|a')_{e_-}}+ \sigma^{(b|b')_{e_-}} \right).$$
\end{ttt}

\begin{proof}
In view of~(\ref{defi_q}) and by symmetry, it is enough to prove the statement of the theorem for $q'$ instead of $q$. Draw geodesics~$\be$ and~$\gamma$ as on Fig.~\ref{vbb'}, where we take~$v:= e_-$.
Similarly to~(\ref{aa'b}) we obtain
$$\big(\beta(j)|\beta(j+1)\big)_{e_-}\ge d(e_-,\beta(j))-1.$$
Using Proposition~\ref{p_bb'} and Lemma~\ref{ab6de},
\begin{eqnarray*}
&&\left|\left<q'[a,b]-q'[a,b_0],e\right>\right|
   = \Big|\Big< \sum_{j= 0}^{d(b_0,b)-1} 
     q'[a,\be(j)]- q'[a,\be(j+1)], e \Big>\Big|\\
&& \le  \sum_{j= 0}^{\infty} 
     \left|\left< q'[a,\be(j)]- q'[a,\be(j+1)], e \right>\right|
   \le  \sum_{j= 0}^{\infty} 
     S_1 \sigma_1^{\big(\be(j)\big|\be(j+1)\big)_{e_-}}\\
&&   \le  \sum_{j= 0}^{\infty} 
     S_1 \sigma_1^{d(e_-,\be(j))-1}
   \le  \sum_{j= 0}^{\infty} 
     S_1 \sigma_1^{d(e_-,\be)+ |j|- 2\de- 1}
   = \frac{S_1\sigma_1^{-2\de-1}}{1-\sigma_1}\,\sigma_1^{d(e_-,\be)}.
\end{eqnarray*}
Similarly,
$$\left|\left<q'[a',b_0]-q'[a',b'],e\right>\right|\le
    \frac{S_1\sigma_1^{-2\de-1}}{1-\sigma_1}\,\sigma_1^{d(e_-,\be)}.$$

Draw a geodesic $\al$ from~$a$ to~$a'$, and let~$a_0\in \al$ be a vertex nearest to~$e_-$.
Interchanging the roles of $a$'s and $b$'s,
a symmetric argument using Proposition~\ref{p_aa'} gives
\begin{eqnarray*}
&& \left|\left<q'[a,b_0]-q'[a_0,b_0],e\right>\right|\le
   \frac{S_2\sigma_2^{-2\de-1}}{1-\sigma_2}\,\sigma_2^{d(e_-,\al)}\mbox{\quad and}\\
&& \left|\left<q'[a_0,b_0]-q'[a',b_0],e\right>\right|\le
   \frac{S_2\sigma_2^{-2\de-1}}{1-\sigma_2}\,\sigma_2^{d(e_-,\al)}.
\end{eqnarray*}
It is an easy exercise using the triangle inequality to see that
$$(a|a')_{e_-}\le d(e_-,\al) \mbox{\qquad and\qquad} (b|b')_{e_-}\le d(e_-,\be).$$
Combining the bounds above, we obtain
\begin{eqnarray*}
&& \left|\left<q'[a,b]-q'[a',b'],e\right>\right|\\
&& \le \left|\left<q'[a,b]-q'[a,b_0],e\right>\right|+
       \left|\left<q'[a,b_0]-q'[a_0,b_0],e\right>\right|+\\
&&     \left|\left<q'[a_0,b_0]-q'[a',b_0],e\right>\right|+
       \left|\left<q'[a',b_0]-q'[a',b'],e\right>\right|\\
&& \le 2 \frac{S_2\sigma_2^{-2\de-1}}{1-\sigma_2}\,\sigma_2^{d(e_-,\al)}+
       2 \frac{S_1\sigma_1^{-2\de-1}}{1-\sigma_1}\,\sigma_1^{d(e_-,\be)}\\
&& \le 2 \frac{S_2\sigma_2^{-2\de-1}}{1-\sigma_2}\,\sigma_2^{(a|a')_{e_-}}+
       2 \frac{S_1\sigma_1^{-2\de-1}}{1-\sigma_1}\,\sigma_1^{(b|b')_{e_-}}.
\end{eqnarray*}
Now the desired inequality follows if we set
$$ \sigma:=\max\{\sigma_1,\sigma_2\} \mbox{\qquad and\qquad}
S:= \max\Big\{  2 \frac{S_1\sigma_1^{-2\de-1}}{1-\sigma_1},\ 
                2\frac{S_2\sigma_2^{-2\de-1}}{1-\sigma_2}
        \Big\}.$$ 
This proves Theorem~\ref{aa'bb'}.
\end{proof}

%%%%%%%%%%%%%%%%%%%%%%%%%%%%%%%%%%%%%%%%%%%%%%%%%%%%%%%%%%%%%%%%%%%%%%
\section{The Ideal Bicombing and a Cocycle}

\subsection{Extending the bicombing}%%%%%%%%%%%%%%%%%%%%%%%%%%%%%%%%%
\label{sec_ideal}%

We fix homological bicombings $q',q$ as above. We proceed now to extend $q$ to infinity.

\medskip

Since $\GG= \G\sqcup \p\G$ is compact, it has a \emph{unique} uniform structure $\mathcal{U}$ compatible with its topology. We shall use the general principle that uniformly continuous maps defined on a dense subset extend continuously if the range is a complete Hausdorff uniform space~\cite[II \S3 \No6]{BourbakiTGI}. For a basepoint $x_0\in\go$, the restriction of $\mathcal{U}$ to $\go$ is given by the base of entourages
$$\big\{(x,y)\in\go\times \go : x=y\ \text{or}\ (x|y)_{x_0}>r\big\},\ \ r>0.$$
Therefore, as this holds for any $x_0$, Theorem~\ref{aa'bb'} states that for all $e\in E$ the function $\left<q[\cdot,\cdot],e\right>$ is uniformly continuous in both variables, and thus uniformly continuous on $\go\times \go$ for the product uniform structure. Recall that any topological vector space has a canonical uniform structure. The preceding comments amount to saying that the map $q:\go\times\go\to \ell^\infty(E)$ is uniformly continuous for the uniform structure corresponding to the topology of pointwise convergence. This topology is generally weaker than the \weak topology, but they coincide on norm-bounded sets (as follows \emph{e.g.} from the Banach\ti Alao\u{g}lu Theorem~\cite[I.3.15]{RudinFA}). By Proposition~\ref{2003de} the map $q'$, hence also $q$, does indeed range in the ball $B$ of radius $2003\,\de^2$. But $B$ is \weak complete and $\go$ is dense in $\go\sqcup\p\G$, so we conclude that $q$ extends (uniquely) to a \weak continuous map
$$q:\ (\go\sqcup\p\G)\times(\go\sqcup\p\G) \lra \ell^\infty(E).$$
We claim that the restriction of $q$ to $\p\G$ is an ideal bicombing as sought for Theorem~\ref{thm_existence}:

\smallskip

Equivariance follows from the uniqueness of the extension. It is easy to check that the image of a closed $|\!|\cdot{|\!|}_1$\ti ball in $\ell^1$ under the inclusion $\ell^1\to\ell^\infty$ is \weak closed in $\ell^\infty$. Therefore, bounded area follows from the fact~(\ref{eq_b_area}) that $q$ has bounded area on $\go$. To be quasigeodesic is a pointwise condition, so it follows from \weak continuity and the corresponding property on $\go$, since we may in particular write $q[\xi,\eta]$ at the \weak limit of $q[a_n,b_n]$ where $a_n,b_n$ tend to $\xi,\eta$ along a fixed geodesic connecting $\xi$ to $\eta$ (unless $\xi=\eta$, in which case $q$ vanishes anyway); indeed with this choice $q[a_n,b_n]$ is supported in a $27\de$\ti neighbourhood of any such geodesic,
by Proposition~\ref{q'}(ii).

Thus it remains only to show that $q$ satisfies the two conditions of Definition~\ref{def_ideal}. For the first condition, fix a vertex $x$ of $\G$; then~(\ref{eq_bic}) shows that $\p q[a,b](x)$ vanishes as $a,b$ tend to the boundary. The condition $\p q[\cdot,\cdot](x)=0$ involves only a fixed finite set of edges, namely the edges incident to $x$. Therefore the \weak continuity of $q$ implies $\p q[\xi,\eta](x)=0$ for all $\xi,\eta\in\p\G$.

\setlength{\unitlength}{1 cm} 
\begin{figure}[h]
  \begin{center}
   \begin{picture}(8,5)
%\graphpaper[1](0,0)(8,5)
% from xi- to xi+
\qbezier(0,0)(2,3)(8,4)
\put(-0.3,-0.3){\footnotesize$\xi_-$}
\put(8.1,4){\footnotesize$\xi_+$}
\qbezier(0,2)(4,3.5)(2,0)
\qbezier(6.5,5)(4,2.9)(7.7,2.7)
\put(2.6,0.6){\footnotesize$V_-$}
\put(7,2.3){\footnotesize$V_+$}
% circle
\put(5.5,3.5){\circle{5}}
\put(5,3.7){\footnotesize$D$}

   \end{picture}
  \end{center}
 \caption{\label{cond(ii)} Illustration for Definition~\ref{def_ideal}~(ii).}
\end{figure}
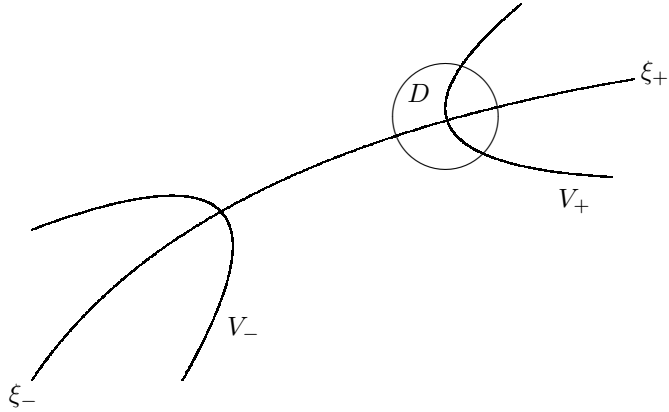

Turning to the second condition, let $\xi_\pm$ and $V_\pm$ be as in Definition~\ref{def_ideal}, see Fig.~\ref{cond(ii)}. Let $D_0$ be the set of vertices $x\in\go$ in the $27\de$\ti neighbourhood of some geodesic from $\xi_-$ to $\xi_+$ such that $x$ is incident to an edge in $V_+$ aswell as to an edge not in $V_+$. Since $V_+$ lies in the complement of a neighbourhood of $\xi_-$ (and since all geodesics connecting $\xi_\pm$ are uniformly close), the set $D_0$ is finite. Let $D\subseteq \G$ be the finite subgraph spanned by $N(D_0,1)$. Let $a^\pm_n$ be sequences tending to $\xi_\pm$ along a given geodesic connecting $\xi_\pm$. Since $q$ is $27\de$\ti quasigeodesic, the condition $\p q[a^-_n, a^+_n]=a^+_n-a^-_n$ shows that
\begin{equation}
\label{eq_supp}%
supp\Big(\p(q[a^-_n, a^+_n]|_{V_+})\Big)\subseteq D\cup\{a^+_n\}\kern1cm\text{for $n$ big enough.}
\end{equation}
By the definition of the boundary operator $\p$, the boundary of any $1$\ti chain (on any subgraph of $\G$) sums to zero. Thus~(\ref{eq_supp}) together with $\p q[a^-_n, a^+_n](a^+_n) = 1$ shows that $\partial (q[\xi_-, \xi_+]|_{V_+})$ is supported on $D$ and sums to $-1$. The argument for $V_-$ is the same. This completes the proof of Theorem~\ref{thm_existence}.\hfill\qedsymbol

\subsection{A cocycle at infinity}%%%%%%%%%%%%%%%%%%%%%%%%%%%%%%

We now construct a non-vanishing cocycle at infinity. The idea is to ``double'' an ideal bicombing into a map $\alpha$ ranging in functions on pairs of edges in such a way that its coboundary $\omega$ never vanishes; this is analogous to a construction given in~\cite{Monod-Shalom1} for trees and $\CAT(-1)$ spaces. 

\begin{ttt}
\label{thm_cocycle}%
Let $\G=(\go,E)$ be a hyperbolic graph of bounded valency. Then there is a \weak continuous $\isom(X)$\ti equivariant alternating cocycle
$$\omega:\ (\p\G)^3 \longrightarrow \ell^1(E\times E)$$
vanishing nowhere on $\partial^3 \G$.
\end{ttt}

\begin{rem}
\label{rem_bounded}%
Such a cocycle will have uniformly bounded $\ell^1$\ti norm on $(\p\G)^3$  since the latter is compact and since \weak compact sets are always bounded in norm~\cite[II\S3]{Day62}. In the construction below, however, the boundedness is granted anyway because we use an ideal bicombing of bounded area.
\end{rem}

First, an independent basic lemma on ideal bicombings which shows that they indeed ``connect'' points at infinity:

\begin{lem}
\label{lemma_non-zero}%
Let $q$ be an ideal quasi-geodesic bicombing. There is a constant $D$ such that for every distinct $\xi_-, \xi_+\in \p\G$, every geodesic $\gamma$ from $\xi_-$ to $\xi_+$ and every vertex $x\in \gamma$ there is an edge $e$ in $\overline{B}(x,D)$
with $\left<q[\xi_-,\xi_+],e\right>\neq 0$.
\end{lem}

\setlength{\unitlength}{1 cm} 
\begin{figure}[h]
  \begin{center}
   \begin{picture}(8,5)
%\graphpaper[1](0,0)(8,5)
% from xi- to xi+
\qbezier(0,1)(4,4)(8,4.2)
\qbezier(0,0.4)(4.5,3.6)(8,3.7)
\qbezier(0,1.6)(4,4.5)(8,4.7)
\put(-0.45,0.75){\footnotesize$\xi_-$}
\put(8.3,4.15){\footnotesize$\xi_+$}
\put(2,2.3){\vector(1,-1){0.377}}
\put(1.8,1.9){\tiny$C$}
% circle
\put(4,3.3){\circle*{0.1}}
\put(4,3.3){\circle{5}}
\put(4,3.3){\vector(1,-1){0.49}}
\put(3.95,2.85){\tiny$D$}
% X+-
\qbezier(4.3,5)(3.3,2)(8,2)
\qbezier(0,4)(5,5)(4,0)
\put(7,2.6){\footnotesize$X_+$}
\put(3,1){\footnotesize$X_-$}

   \end{picture}
  \end{center}
 \caption{Illustration for Lemma~\ref{lemma_non-zero}.}
\end{figure}
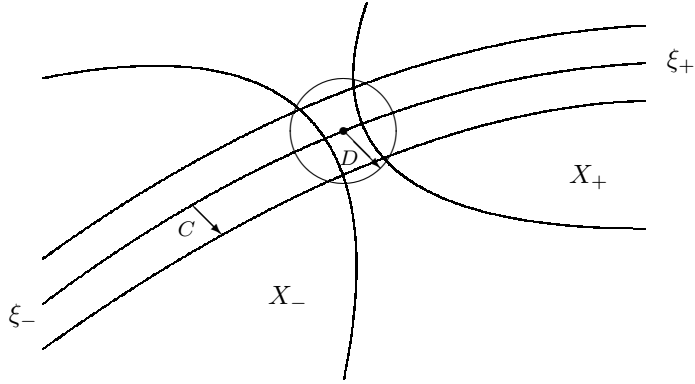

\begin{proof}[Proof of the lemma]
Let $C$ be a constant as in the definition of quasi-geodesic ideal bicombings and fix $D\geq C+2+3\de$. Given a geodesic $\gamma: {\bf R}\to \G$ and $t\in{\bf R}$, let $X^+\se E$ (respectively $X^-$) be the set of edges $e\in E$ such that some $s\geq t+D-C-1$ (resp. some $s\leq t-(D-C-1)$) realizes the minimum over $\bf{R}$ of the function $s\mapsto d(\gamma(s), e)$. Then $X^+\cup X^-\cup \overline{B}(\gamma(t),D)$ covers $N(\gamma,C)$. On the other hand, by $\de$\ti hyperbolicity, $X^+\cap X^- = \varnothing$. Apply this to any geodesic as in the statement with $x=\gamma(t)$; assuming for a contradiction $q[\xi_-,\xi_+]|_{\overline{B}(x,D)}=0$, we deduce $\partial(q[\xi_-, \xi_+]|_{X^\pm})=0$ since $0=\partial(q[\xi_-, \xi_+]) = \partial(q[\xi_-, \xi_+]|_{X^-}) + \partial(q[\xi_-, \xi_+]|_{X^+})$  and $X^+\cap X^- = \varnothing$. This contradicts the definition of an ideal bicombing because it can be readily checked that $X^\pm$ is a graph neighbourhood of~$\xi_\pm$ using $\de$\ti hyperbolicity.
\end{proof}

\begin{proof}[Proof of Theorem~\ref{thm_cocycle}]
Let $q$ be an equivariant \weak continuous quasi-geodesic ideal bicombing of bounded area, as granted by Theorem~\ref{thm_existence}, and let $C$ be as above. We may assume $q$ alternating upon replacing $q[\xi,\eta]$ with $(q[\xi,\eta] - q[\eta,\xi])/2$ (the particular bicombing $q$ constructed above is alternating anyway). Given a constant $R$ we define $\alpha:(\p\G)^2 \to \ell^\infty(E\times E)$ by
$$\alpha(\xi,\eta)(e,e') =\begin{cases}
\left<q[\xi,\eta],e\right> \left<q[\xi,\eta],e'\right> &\text{if $d(e,e')\leq R$,}\\
0 &\text{otherwise.}
\end{cases}$$
Now we set $\omega = d\alpha$, that is, by alternation, $\omega(\xi,\eta,\zeta):= \alpha(\xi,\eta)+\alpha(\eta,\zeta)+\alpha(\zeta,\xi)$. We shall check (i)~that $\omega$ ranges in $\ell^1(E\times E)$ and even in a closed ball thereof, and (ii)~that (for a suitable choice of $R$) it vanishes nowhere on $\partial^3 \G$, as the other properties follow from the definition. In particular, \weak continuity of $\omega$ for $\ell^1$ (in duality to $\co$) follows from pointwise continuity because it ranges within a closed $|\!|\cdot{|\!|}_1$\ti ball (Remark~\ref{rem_bounded} and~(i) below); we used this trick for $\ell^\infty$ in Section~\ref{sec_ideal}.

%\setlength{\unitlength}{1 cm} 
%\begin{figure}[h]
%  \begin{center}
%   \begin{picture}(8,4)
%\graphpaper[1](0,0)(8,4)
% from xi to eta
%\qbezier(0,2)(3,2.05)(4.3,2.4)
%\qbezier(4.3,2.4)(5.5,2.7)(8,4)
%\put(-0.2,2){\footnotesize$\xi$}
% from xi to zeta
%\qbezier(0,2)(3,1.95)(4.3,1.6)
%\qbezier(4.3,1.6)(5.5,1.3)(8,0)
%\put(8,-.2){\footnotesize$\zeta$}
% from eta to zeta
%\qbezier(8,4)(5,2.5)(5,2)
%\qbezier(5,2)(5,1.5)(8,0)
%\put(8.1,4){\footnotesize$\eta$}
% circle
%\put(4.5,2){\circle{1.5}}
%\put(4.5,2){\vector(-1,2){0.31}}
%\put(4.45,2.2){\footnotesize$r$}
% A_xi
%\qbezier(0.5,3)(4.2,3)(4.2,2)
%\qbezier(4.2,2)(4.2,1)(0.5,1)
%\put(2,3.1){\footnotesize$A_\xi$}
% A_eta
%\qbezier(6,4)(4,3)(4.8,2.3)
%\qbezier(4.8,2.3)(5.6,1.6)(8,3)
%\put(5,3.9){\footnotesize$A_\eta$}
% A_zeta
%\qbezier(6,0)(4,1)(4.8,1.7)
%\qbezier(4.8,1.7)(5.6,2.4)(8,1)
%\put(5,0.1){\footnotesize$A_\zeta$}

%   \end{picture}
%  \end{center}
% \caption{\label{point(i)} Point~(i).}
%\end{figure}

\textbf{(i)}~By hyperbolicity, there is a constant $r$ such that for any distinct $\xi,\eta,\zeta\in\p\G$ and any three geodesics joining them, the union of the $C$\ti neighbourhoods of the three geodesics can be covered by some ball of radius $r$ together with three subgraphs $A_\xi$, $A_\eta$, $A_\zeta$ such that $A_\xi$ intersects only the $C$\ti neighbourhoods of the two geodesics leading to $\xi$~--- and likewise for $A_\eta$, $A_\zeta$. Thus, by the symmetry of the situation, it is enough to show that the restriction of $\omega(\xi,\eta,\zeta)$ to pairs of edges in $A_\xi$ is summable (observe indeed that there are only uniformly finitely many pairs $(e,e')$ with $\omega(\xi,\eta,\zeta)(e,e')\neq 0$ for which $e$ is in $A_\xi$ but $e'$ is not). Let $M$ be a bound on the number of edges in any ball of radius $R+1$; we have:
\begin{multline*}
|\!|\omega(\xi,\eta,\zeta)|_{A_\xi}{|\!|}_1 = \sum_{e,e' \in A_\xi} 
  \big|\alpha(\xi,\eta)(e,e') - \alpha(\xi,\zeta)(e,e') \big| \leq\\
\sum_{e \in A_\xi}\sum_{e' \in N(e,R)}
 \big|\left<q[\xi,\eta],e\right>\big|\cdot 
 \big|\left<q[\xi,\eta],e'\right> - \left<q[\xi,\zeta],e'\right>\big| +\\
\ \kern1cm  \sum_{e' \in A_\xi}\sum_{e \in N(e',R)}
   \big|\left<q[\xi,\eta],e\right> - \left<q[\xi,\zeta],e\right>\big|\cdot 
   \big|\left<q[\xi,\zeta],e'\right>\big| \leq\\
\ \kern1cm 2 M |\!|q{|\!|}_\infty \sum_{e \in A_\xi}
   \big|\left<q[\xi,\eta],e\right> - \left<q[\xi,\zeta],e\right>\big| = \\
2 M |\!|q{|\!|}_\infty \sum_{e \in A_\xi}
  \big|\left<q[\xi,\eta],e\right> +
  \left<q[\eta,\zeta],e\right> + 
  \left<q[\zeta,\xi],e\right>\big|.
\end{multline*}
The last term is finite and uniformly bounded since $q$ has bounded area.

\setlength{\unitlength}{1 cm} 
\begin{figure}[h]
  \begin{center}
   \begin{picture}(10,5)
%\graphpaper[1](0,0)(10,5)
% from xi to eta
\qbezier(0,1)(3,0.92)(4.3,1.33)
\qbezier(4.3,1.33)(5.3,1.6)(10,5)
\put(5,2.1){\footnotesize$\gamma_{\xi\eta}$}
\put(-0.2,1){\footnotesize$\xi$}
% from xi to zeta
\qbezier(0,1)(4,1)(8,0)
\qbezier(0,1.5)(4,1.5)(8,0.5)
\put(7.8,0.05){\vector(0,1){0.5}}
\put(7.4,0.25){\footnotesize$C$}
\put(4,0.4){\footnotesize$\gamma_{\zeta\xi}$}
\put(8.1,-0.1){\footnotesize$\zeta$}
% from eta to zeta
\qbezier(10,5)(5,1.3)(5,1)
\qbezier(5,1)(5,0.7)(8,0)
\put(6.4,2){\footnotesize$\gamma_{\eta\zeta}$}
\put(10.1,5){\footnotesize$\eta$}
% circle
\put(7.97,3.55){\circle{5}}
\put(7.97,3.55){\circle*{0.1}}
\put(7.97,3.55){\vector(-1,0){.7}}
\put(7.5,3.65){\footnotesize$D$}
\put(8.25,3.2){\footnotesize$e$}
\qbezier(8.25,3.1)(8.175,3.125)(8.1,3.15)
   \end{picture}
  \end{center}
 \caption{\label{point(ii)} Point~(ii).}
\end{figure}
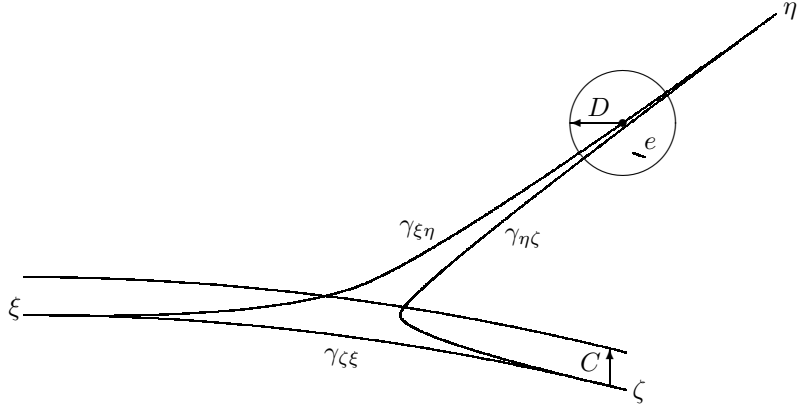

\textbf{(ii)}~Let $D$ be as in Lemma~\ref{lemma_non-zero} and choose $L>2(C+D+\de)$. We claim that whenever $R>2L+2D$ the cocycle $\omega$ vanishes nowhere on $\partial^3 \G$. Indeed, choose any three distinct points $\xi,\eta,\zeta$ in $\p\G$. Let $\gamma_{\xi\eta}$ be a geodesic from $\xi$ to $\eta$, and similarly $\gamma_{\eta\zeta}, \gamma_{\zeta\xi}$
(see Fig.~\ref{point(ii)}). 
There is $t\in{\bf R}$ such that $\gamma_{\xi\eta}(t)$ is at distance at most $\de$ from $\gamma_{\eta\zeta}$ aswell as from $\gamma_{\zeta\xi}$. By Lemma~\ref{lemma_non-zero}, we can find an edge $e$ in the ball $\overline{B}(\gamma_{\xi\eta}(t+L),D)$ with $\left<q[\xi,\eta],e\right>\neq 0$. Likewise, there is an edge $e'$ in $\overline{B}(\gamma_{\xi\eta}(t-L),D)$ with $\left<q[\xi,\eta],e'\right>\neq 0$. However, by the choice of $t$ and $L$ we have $\overline{B}(\gamma_{\xi\eta}(t+L),D)\cap N(\gamma_{\zeta\xi},C) =\varnothing$ and thus $\left<q[\zeta,\xi],e\right>=0$. Likewise, $\left<q[\eta,\zeta],e'\right>=0$. Therefore, $\omega(\xi,\eta,\zeta)(e,e') = \left< q[\xi,\eta],e\right> \left<q[\xi,\eta],e'\right>$ is non-zero because $d(e,e')<R$.
\end{proof}

\section{Superrigidity and remaining proofs}%%%%%%%%%%%%%%%%%%%%%%%%%%%%%%

One of the reasons why one has to work on the boundary at infinity is that there is no known technique to show that general cocycles defined on the group itself yield non-trivial classes in bounded cohomology when one is dealing with arbitrary coefficients. Extending cocycles to infinity settles this problem, see~\ref{sec_poisson} below.

%%%%%%%%%%%%%%%%%%%%%%%%%%%%%%%%%%%%%%%%%%%%%%%%%%%%%%%%%%%%%%%%%%%%%%%%%%%%%%%%%%%%%%%%%%%%%%%%%%%%
\subsection{More general spaces}
\label{sec_general_spaces}%

In order to tackle more general hyperbolic spaces, we establish the following alternative:

\begin{ttt}
\label{thm_struct_hyp}%
Let $X$ be a proper geodesic hyperbolic space and $H<\isom(X)$ a non-elementary closed subgroup acting cocompactly on $X$. Then either:

\smallskip

\noindent(i) $H$ has a proper non-elementary vertex-transitive action on a hyperbolic graph of bounded valency (thus quasi-isometric to $X$), or

\smallskip

\noindent(ii) there is a finite index open subgroup $H^*\lhd H$ and a compact normal subgroup $K\lhd H$ contained in $H^*$ such that $H^*/K$ is (isomorphic to) a connected simple Lie group of rank one.
\end{ttt}

We first give the

\begin{proof}[Proof of Proposition~\ref{prop_def_elem}]
Assume that $H$ fixes a probability measure $\mu$ on $\p X$. If the support of $\mu$ contains at least three points, then $H$ is relatively compact~\cite[5.3]{Adams96} and hence preserves a compact set (any orbit of the closure of $H$) in $X$. Otherwise, the support of $\mu$ is an $H$\ti invariant set of one or two points in $\p X$. The converse is clear as the stabiliser of a compact set is compact and hence fixes a measure on $\p X$.

Condition~(iii) always implies~(i), and the converse follows from universal amenability of the boundary action, which is granted by a result of Adams~\cite[6.8]{Adams96} when $X$ has at most exponential growth. Note that the growth condition is automatically satisfied for graphs of bounded valency, and in the case of a proper geodesic space with cocompact $\isom(X)$\ti action it can be deduced \emph{e.g.} from Lemma~\ref{lemma_comp_gen} below.
\end{proof}

For lack of a reference, we prove the following:

\begin{lem}
\label{lemma_comp_gen}%
Let $X$ be a proper geodesic metric space and $H<\isom(X)$ a closed subgroup. If $H$ acts cocompactly on $X$, then $H$ is compactly generated.
\end{lem}

\begin{proof}[Proof of the lemma]
Fix a point $x\in X$. There is $r>0$ such that $H B(x,r)=X$. It follows, since $\overline{B}(x,2r)$ is compact, that there is a finite set $L\subseteq H$ such that $\overline{B}(x,2r)\subseteq L B(x,r)$. Let $C\subseteq H$ be the compact set of elements $h\in H$ such that $hx\in \overline{B}(x,r)$. We claim that $H$ is generated by $C\cup L$. Indeed, the fact that $X$ is geodesic implies by induction on $n$ that $B(x, nr)\subseteq <L>\, B(x,r)$. Thus $<L> B(x,r) = X$ and the claim follows.
\end{proof}

\begin{proof}[Proof of Theorem~\ref{thm_struct_hyp}]
Let $K\lhd H$ be the maximal amenable normal subgroup; we claim that $K$ is compact. Indeed, $K$ fixes a probability measure on $\p X$. If the support $F$ of this measure contains at least three points, then $K$ is compact~\cite[5.3]{Adams96}. If on the other hand $F$ contains one or two points, then $H$ cannot preserve the set $F$ since $H$ is non-elementary. Thus there is $h\in H$ such that $hF\cup F$ contains at least three points. But since $h$ normalises $K$, the set $hF\cup F$ is preserved by $K$. This implies again that $K$ is compact for the same reason since $K$ fixes the uniform distribution on that set.

Let now $L=H/K$ and let $L^0$ be its connected component. As explained \emph{e.g.} in~\cite[11.3.4]{Monod}, the solution to Hilbert's fifth problem and the classification of outer automorphisms of adjoint semi-simple Lie groups implies that there is a finite index (normal) open subgroup $H^*\lhd H$ containing $K$ such that $H^*/K$ splits as a direct product $L^0\times Q$ where $Q$ is totally disconnected and $L^0$ is a semi-simple connected Lie group without compact factors. Assume first that $L^0$ is non-compact; we claim that $Q$ is trivial and $L^0$ simple.

Indeed, otherwise $H^*/K$ would be a direct product of two non-compact groups since $H^*/K$ has no normal compact subgroups (see~\cite[11.3.3]{Monod}). Let $H_1, H_2 < H^*$ be the preimages of these two factors. Since $H_1$ is non-compact, its limit set $S\subseteq \p X$ is non-empty. But $H_2$ acts trivially on $S$ because the $H_i$ commute up to the compact subgroup $K$; therefore, $S$ contains only one or two points since otherwise $H_2$ would be compact. Both $H_i$ preserve $S$, so that now $H^*$ would be elementary; in view of the first caracterisation in Proposition~\ref{prop_def_elem}, this implies that $H$ is elementary, a contradiction.

Thus, in the case $L^0$ non-compact, it remains only to see that $H^*/K$ has rank one; this is the case because the associated symmetric space is quasi-isometric to $X$.

The other case is when $L^0$ is compact, and hence trivial. Thus $L=H/K$ is totally disconnected; since it is compactly generated by Lemma~\ref{lemma_comp_gen}, it admits a locally finite Schreier graph (see \emph{e.g.}~\cite[p.~150]{Monod}). The $H$\ti action on this graph is proper and vertex-transitive, hence the graph has bounded valency and is quasi-isometric to $X$, which implies both hyperbolicity and non-elementarity.
\end{proof}

We are going to need the existence of a distance-like function for pairs of pairs of points at infinity. Since the invariant cross-ratio construction (as \emph{e.g.} in~\cite{Paulin96}) is not continuous even for hyperbolic groups, we construct a continuous invariant analogue by applying the following lemma to $C=\p X$ and $G<\isom(X)$:

\begin{lem}
\label{lemma_cross}%
Let $G$ be a locally compact group with a continuous action on a compact second countable space $C$ such that the diagonal action on the space of distinct triples of $C$ is proper. Denote by $D_2$ the space of subsets of cardinality two in $C$ and by $D$ the space of distinct (not necessarily disjoint) pairs $(a,b)\in D_2\times D_2$. Then there is a continuous $G$\ti invariant function $\delta:D\to{\bf R}_+$ such that $\delta(a,b)=0$ if and only if $a\cap b\neq\varnothing$.
\end{lem}

\begin{proof}
The space $C$ is metrisable, so let $\ro$ be a metric inducing the topology. For two disjoint elements $a=\{\alpha_1, \alpha_2\}$ and $b=\{\beta_1, \beta_2\}$ of $D_2$, let
$$\delta'(a,b) = \left|\ln\frac{\ro(\alpha_1, \beta_1)\cdot\ro(\alpha_2, \beta_2)}{\ro(\alpha_1, \beta_2)\cdot\ro(\alpha_2, \beta_1)}\right|^{-1}.$$
This extends to a continuous function $\delta':D\to{\bf R}_+$ such that $\delta'(a,b)=0$ if and only if $a\cap b\neq\varnothing$; indeed the definition of $D_2$ excludes the case where both numerator and denominator are zero. Our assumption implies that $G$ acts properly on $D$. Thus, fixing a left Haar measure on $G$, there exists a generalized Bruhat function, that is, a continuous function $h:D\to {\bf R}_+$ such that (i)~for all $(a,b)\in D$ we have $\int_G h(g^{-1}(a,b))\,dg = 1$ and (ii)~for every compact subset $L\subset D$ the intersection of the support of $h$ with the saturation $GL$ is compact (see \emph{e.g.}~\cite[4.5.4]{Monod}). The function $\delta:D\to {\bf R}_+$ defined by
$$\delta(a,b) = \int_G h(g^{-1}(a,b))\delta'(g^{-1}(a,b))\,dg$$
is $G$\ti invariant and continuous (for continuity, see \emph{e.g.} the proof of~\cite[4.5.5]{Monod}; the integrand vanishes outside a compact set). Moreover, $\delta(a,b)>0$ if $a,b$ are disjoint because of~(ii).
\end{proof}

\subsection{Boundary theory and bounded cohomology}%%%%%%%%%%%%%%%%%%%%%%%%%%%%%%
\label{sec_poisson}%

We shall prove the following more general form of Theorem~\ref{thm_coho}. For background on (continuous) bounded cohomology, see~\cite{Burger-Monod3,Monod}.

\begin{ttt}
\label{thm_coho_bis}%
Let $\G$ be a hyperbolic graph of bounded valency and $H<\isom(\G)$ a non-elementary closed subgroup. Then $\hbc^2(H,L^p(H))$ is non-zero for all $1\leq p < \infty$.
\end{ttt}

Using Theorem~\ref{thm_struct_hyp}, we can combine Theorem~\ref{thm_coho_bis} with~\cite{Monod-Shalom1} and deduce:

\begin{ccc}
\label{cor_coho_bis}%
Let $X$ be a hyperbolic proper geodesic metric space such that $\isom(X)$ acts cocompactly. For any non-elementary closed subgroup $G<\isom(X)$, the space $\hbc^2(G,L^2(G))$ is non-zero.
\end{ccc}

\begin{rem}
In fact the proof works for coefficients in $L^p(G)$ for all $1<p<\infty$, but contrary to Theorem~\ref{thm_coho_bis} the method does not apply to $p=1$.
\end{rem}

\begin{proof}[Proof of the corollary]
Let $H=\isom(X)$. Since $H$ contains the non-elementary subgroup $G$, it is itself non-elementary and we may apply Theorem~\ref{thm_struct_hyp}. In case~(i), we are done by Theorem~\ref{thm_coho_bis} applied to $G$. In case~(ii), there is a finite index subgroup of $G$ with a proper non-elementary action on a rank one symmetric space, a situation covered by~\cite{Monod-Shalom1}.
\end{proof}

We need the notion of (double) $\mathfrak{X}^\mathrm{sep}$\ti ergodicity introduced by Burger-Monod in~\cite{Burger-Monod3}. Recall that if $H$ is any locally compact $\sigma$\ti compact group acting on a measure space $(B,\nu)$ by measurable transformations leaving $\nu$ quasi-invariant, then $B$ is called {\sf $\mathfrak{X}^\mathrm{sep}$\ti ergodic} if for every separable coefficient $H$\ti module $F$, any measurable $H$\ti equivariant function $f:B\to F$ is essentially constant. (Here, a {\sf coefficient module} is the contragredient of a continuous isometric representation on a separable Banach space.) Further, $B$ is said {\sf doubly} $\mathfrak{X}^\mathrm{sep}$\ti ergodic if the diagonal action on $B\times B$ is $\mathfrak{X}^\mathrm{sep}$\ti ergodic. It is proved in~\cite{Burger-Monod3,Monod} that if $A$ is compactly generated then a certain Poisson boundary is a doubly ergodic $H$\ti space. This has been generalized by Kaimanovich~\cite{Kaimanovich03} to the $\sigma$\ti compact case. On the other hand, Poisson boundaries are amenable in the sense of Zimmer, see~\cite{Zimmer78b}. (For comments about amenable spaces for possibly not second countable groups, see~\cite[5.3]{Monod}; as to measurability issues, all our $H$\ti actions considered here factor through a second countable quotient.) Using the cohomological characterization of amenability established in~\cite{Burger-Monod3,Monod}, it follows that if $B$ is a Poisson boundary as above, then for any separable coefficient $H$\ti module $F$ there is a canonical isometric isomorphism
\begin{equation}
\label{eq_h2}%
\hbc^2(H,F)\cong ZL^\infty_\mathrm{alt}(B^3,F)^H,
\end{equation}
\noindent
where the right hand side denotes the space of equivariant essentially bounded measurable alternating cocycles; see~\cite{Burger-Monod3,Monod}.

\begin{ppp}
\label{prop_Furst}%
Let $X$ be a proper hyperbolic geodesic metric space and $H$ a locally compact $\sigma$\ti compact group with a continuous non-elementary action on $X$. Let $B$ be a doubly ergodic amenable $H$\ti space. Then (after possibly discarding a null-set in $B$) there is a Borel $H$\ti equivariant map $f:B\to \p X$.
\end{ppp}

\begin{proof}
An examination of the proof given in~\cite{Monod-Shalom1} for the $\CAT(-1)$ case (in the more general cocycle setting) shows that the only two ingredients about $\p X$ that are really needed are (i)~the existence of a function~$\delta$ as in Lemma~\ref{lemma_cross} (for $C=\p X$ and $G$ the closure in $\isom(X)$ of the image of $H$), and (ii)~the fact~\cite[5.3]{Adams96} that $\isom(X)$ acts properly on the space of probabilities on $\p X$ whose support contains at least three points.
\end{proof}

\begin{proof}[End of the proof of Theorem~\ref{thm_coho_bis}]
Let $(B,\nu)$ be a doubly $\mathfrak{X}^\mathrm{sep}$\ti ergodic amenable $H$\ti space and $f:B\to \p X$ as in Proposition~\ref{prop_Furst}. Let now $F=\ell^1(E\times E)$ and $\omega: (\p\G)^3\to F$ be as in Theorem~\ref{thm_cocycle}; consider the resulting map $f^*\omega:B^3\to F$. We observe that $F$ is a separable coefficient module. The support of $f_*\nu$ is an $H$\ti invariant set in $\p\G$, so that it contains at least three points since $H$ is non-elementary. The Theorem of Fubini-Lebesgue implies that $\partial^3\G$ is not $(f_*\nu)^3$\ti null, hence $f^*\omega$ defines a non-zero element in the right hand side of~(\ref{eq_h2}) because $\omega$ is nowhere vanishing on $\p^3\G$. We claim that $\hbc^2(H,\ell^1(E))$ is also non-vanishing. Indeed this follows readily from~(\ref{eq_h2}) because $\ell^1(E\times E)$ is isomorphic (as coefficient module) to the completed sum of copies of $\ell^1(E)$. Similarily, if $K<H$ is the stabiliser of a given edge in $E$, the space $\hbc^2(H,\ell^1(H/K))$ is non-zero. Using~\cite[11.4.1]{Monod} and the fact that $\ell^p(H/K)$ is dual even for $p=1$ since $K$ is open, this implies
\begin{equation}
\label{eq_lp_nonzero}%
\hbc^2\big(H,\ell^p(H/K)\big)\neq 0\ \ \ \forall\,1\leq p<\infty.
\end{equation}
(This is false for $p=\infty$, and indeed~\cite[11.4.1]{Monod} does not apply by lack of separability.) Since $K$ is compact, the submodule $\ell^p(H/K)$ is equivariantly complemented in $L^p(H)$. Now the conclusion follows from~\cite[8.2.9]{Monod}.
\end{proof}

Observe that the properties of the map of Proposition~\ref{prop_Furst} are the only point of the whole proof where we use non-elementarity. If on the other hand $H$ is elementary, then it is amenable by Proposition~\ref{prop_def_elem} and thus $\hbc^2(H,\ell^p(E))$ vanishes altogether since $\ell^p(E)$ is a dual module (compare \emph{e.g.}~\cite[2.5]{Johnson}).

\begin{proof}[End of the proof of Theorem~\ref{thm_coho}]
The situation of Theorem~\ref{thm_coho} corresponds to a homomorphism $\pi:\Gamma\to \isom(\G)$ with $\Gamma_0:=\pi(\Gamma)$ a discrete non-elementary subgroup. Thus Theorem~\ref{thm_coho_bis} implies that $\hb^2(\Gamma_0, \ell^p(\Gamma_0))$ is non-zero. The conclusion follows because the kernel of $\pi$ is finite.
\end{proof}

\begin{proof}[Proof of Corollary~\ref{cor_general_hyp}]
Keep the notation of Corollary~\ref{cor_general_hyp} and let $H=\isom(X)$; we may apply Theorem~\ref{thm_struct_hyp} since $H$ contains a non-elementary subgroup, namely the image of $\Gamma$. In case~(i), we are done by Theorem~\ref{thm_coho}. In case~(ii), there is a finite index subgroup of $\Gamma$ with a proper non-elementary action on a rank one symmetric space. Then the results of~\cite{Monod-Shalom1} imply that $\Gamma$ is in~$\mathcal{C}_\mathrm{reg}$.
\end{proof}

\begin{rem}
\label{rem_Ilya}%
We mentioned in the introduction that there are uncountably many non-isomorphic subgroups of hyperbolic groups. Here is the sketch of the argument given to us by I.~Kapovich and P.~Schupp; we thank them for communicating it. Let $A=\bigoplus_p {\bf Z}/p$ be the sum of all cyclic groups of prime order. Since $A$ is countable and recursively presented, it can be embedded in a finitely presented group $Q$ by a result of Higman's~\cite{Higman}. The construction of Rips in~\cite{Rips82} can be modified so as to give a hyperbolic group $\Gamma$ with a homomorphism $\pi$ onto $Q$ such that the kernel $K$ of $\pi$ is perfect.
$$\xymatrix{
A \ar@{^(->}[r]^i& Q & \Gamma\rhd K \ar@{->>}[l]_{\pi}
}$$
Now observe that for any subgroup $B<A$ the subgroup $\Gamma_B=\pi^{-1}(i(B))$ of $\Gamma$ has Abelianisation isomorphic to $B$ since $K$ is perfect. Since there is a continuum of non-isomorphic $B<A$ (by taking sums restricted to any set of primes), the conclusion follows.
\end{rem}

\begin{proof}[Proof of Theorem~\ref{thm_superrigidity}]
Keep the notation of Theorem~\ref{thm_superrigidity} and suppose $H$ non-amenable, so that it is non-elementary (see the Proof of Proposition~\ref{prop_def_elem} above). Assume first that we are in case~(i) with $X=\G=(\go,E)$. Arguing as in the proof of Theorem~\ref{thm_coho_bis} but with a doubly $\mathfrak{X}^\mathrm{sep}$\ti ergodic amenable $\Lambda$\ti space (and $\Lambda$ replacing $H$), we deduce that $\hb^2(\Lambda, \ell^2(E))$ is non-zero, and thus $\hb^2(\Lambda, L^2(H))$ does not vanish either. In case~(ii), Theorem~\ref{thm_struct_hyp} allows us to apply results from~\cite{Monod-Shalom1} and deduce again that $\hb^2(\Lambda, L^2(H))$ is non-zero. Now the superrigidity formula for bounded cohomology of irreducible lattices (Theorem~16 in~\cite{Burger-Monod3}) implies that there is a non-zero $\Lambda$\ti invariant subspace of $L^2(H)$ on which the $\Lambda$\ti representation extends to a continuous $G$\ti representation that factors through $G\to G_i$ for some $i$. (In order to relax the compact generation assumption in~\cite{Burger-Monod3} to $\sigma$\ti compactness, apply~\cite{Kaimanovich03}.) Using the arguments of the proof of Theorem~0.3 in~\cite{Shalom00}, one deduces that there is a compact normal subgroup $K\lhd H$ such that the homomorphism $\Gamma\to H/K$ extends to a continuous homorphism $G\to H/K$ factoring through $G_i$.
\end{proof}

%%%%%%%%%%%%%%%%%%%%%%%%%%%%%%%%%%%%%%%%%%%%%%%%%%%%%%%%%%%%%%%%%%%%%%
\ifx\undefined\bysame
\newcommand{\bysame}{\leavevmode\hbox to3em{\hrulefill}\,}
\fi

\smallskip

\noindent
I. M.: Department of Mathematics, University of Illinois at Urbana-Champaign, Urbana, IL 61801 USA, {\tt mineyev@math.uiuc.edu}\\

\noindent
N. M.: Department of Mathematics, University of Chicago, Chicago, IL 60637 USA, \smash{\tt monod@uchicago.edu}\\

\noindent
Y. S.: School of Mathematical Sciences, Tel-Aviv University, Tel-Aviv 69978, Israel, {\tt yeshalom@post.tau.ac.il}

\begin{thebibliography}{GdlH90}

\bibitem[Ada94]{Adams94a}
Scot Adams, {\em Indecomposability of equivalence relations generated by word
  hyperbolic groups}, Topology {\bf 33} (1994), no.~4, 785--798.

\bibitem[Ada95]{Adams95}
Scot Adams, {\em {Some new rigidity results for stable orbit equivalence}},
  Ergodic Theory Dyn. Syst. {\bf 15} (1995), no.~2, 209--219.

\bibitem[Ada96]{Adams96}
Scot Adams, {\em {Reduction of cocycles with hyperbolic targets}}, Ergodic
  Theory Dyn. Syst. {\bf 16} (1996), no.~6, 1111--1145.

\bibitem[BF02]{Bestvina-Fujiwara}
Mladen Bestvina and Koji Fujiwara, {\em Bounded cohomology of subgroups of
  mapping class groups}, Geom. Topol. {\bf 6} (2002), 69--89 (electronic).

\bibitem[BM02]{Burger-Monod3}
Marc Burger and Nicolas Monod, {\em Continuous bounded cohomology and
  applications to rigidity theory (with an appendix by {M}.~{B}urger and
  {A}.~{I}ozzi)}, Geom. Funct. Anal. {\bf 12} (2002), no.~2, 219--280.

\bibitem[Bou71]{BourbakiTGI}
Nicolas Bourbaki, {\em \'{E}l\'ements de math\'ematique. {T}opologie
  g\'en\'erale. {C}hapitres 1 \`a 4}, Hermann, Paris, 1971.

\bibitem[Cha94]{Champetier}
Christophe Champetier, {\em Petite simplification dans les groupes
  hyperboliques}, Ann. Fac. Sci. Toulouse Math. (6) {\bf 3} (1994), no.~2,
  161--221.

\bibitem[Day62]{Day62}
Mahlon~M. Day, {\em Normed linear spaces}, Academic Press Inc., Publishers, New
  York, 1962.

\bibitem[EF97]{Eskin-Farb1}
Alex Eskin and Benson Farb, {\em Quasi-flats and rigidity in higher rank
  symmetric spaces}, J. Amer. Math. Soc. {\bf 10} (1997), no.~3, 653--692.

\bibitem[EF98]{Eskin-Farb2}
Alex Eskin and Benson Farb, {\em Quasi-flats in ${H}\sp 2\times {H}\sp 2$}, Lie
  groups and ergodic theory (Mumbai, 1996) (Bombay), Tata Inst. Fund. Res.,
  Bombay, 1998, pp.~75--103.

\bibitem[FLM01]{Farb-Lubotzky-Minsky}
Benson Farb, Alexander Lubotzky, and Yair~N. Minsky, {\em Rank-1 phenomena for
  mapping class groups}, Duke Math. J. {\bf 106} (2001), no.~3, 581--597.

\bibitem[Fuj98]{Fujiwara98}
Koji Fujiwara, {\em The second bounded cohomology of a group acting on a
  {G}romov-hyperbolic space}, Proc. London Math. Soc. (3) {\bf 76} (1998),
  no.~1, 70--94.

\bibitem[Gab02]{GaboriauL2}
Damien Gaboriau, {\em Invariants {$l\sp 2$} de relations d'\'equivalence et de
  groupes}, Publ. Math. Inst. Hautes \'Etudes Sci. (2002), no.~95, 93--150.

\bibitem[GdlH90]{Ghys-Harpe}
{\'E}tienne Ghys and Pierre de~la Harpe (eds.), {\em Sur les groupes
  hyperboliques d'apr\`es {M}ikhael {G}romov}, Birkh\"auser Verlag, Basel,
  1990, Papers from the Swiss Seminar on Hyperbolic Groups held in Bern, 1988.

\bibitem[Gro87]{Gromov87}
Mikha{\"\i}l Gromov, {\em Hyperbolic groups}, Essays in group theory (New
  York), Math. Sci. Res. Inst. Publ., vol.~8, Springer, New York, 1987,
  pp.~75--263.

\bibitem[Gro93]{Gromov91}
Mikha{\"\i}l Gromov, {\em Asymptotic invariants of infinite groups}, Geometric
  group theory, Vol.\ 2 (Sussex, 1991) (Cambridge), Cambridge Univ. Press,
  Cambridge, 1993, pp.~1--295.

\bibitem[Hig61]{Higman}
Graham Higman, {\em Subgroups of finitely presented groups}, Proc. Roy. Soc.
  Ser. A {\bf 262} (1961), 455--475.

\bibitem[Joh72]{Johnson}
Barry~E. Johnson, {\em {Cohomology in Banach algebras}}, Mem. Am. Math. Soc.
  {\bf 127} (1972).

\bibitem[Kai]{Kaimanovich03}
Vadim~A. Kaimanovich, {\em Double ergodicity of the {P}oisson boundary and
  applications to bounded cohomology}, Preprint to appear in GAFA.

\bibitem[KH]{Kechris-Hjorth}
Alexander~S. Kechris and Greg Hjorth, {\em Rigidity theorems for actions of
  product groups and countable borel equivalence relations}, Preprint.

\bibitem[KL97]{Kleiner-Leeb}
Bruce Kleiner and Bernhard Leeb, {\em Rigidity of quasi-isometries for
  symmetric spaces and {E}uclidean buildings}, Inst. Hautes \'Etudes Sci. Publ.
  Math. (1997), no.~86, 115--197 (1998).

\bibitem[Min01]{Mineyev01}
Igor Mineyev, {\em Straightening and bounded cohomology of hyperbolic groups},
  Geom. Funct. Anal. {\bf 11} (2001), no.~4, 807--839.

\bibitem[MM99]{Masur-Minsky}
Howard~A. Masur and Yair~N. Minsky, {\em Geometry of the complex of curves.
  {I}. {H}yperbolicity}, Invent. Math. {\bf 138} (1999), no.~1, 103--149.

\bibitem[Mon01]{Monod}
Nicolas Monod, {\em {Continuous bounded cohomology of locally compact groups}},
  {Lecture Notes in Mathematics 1758, Springer, Berlin}, 2001.

\bibitem[MSa]{Monod-Shalom1}
Nicolas Monod and Yehuda Shalom, {\em Cocycle superrigidity and bounded
  cohomology for negatively curved spaces}, Preprint.

\bibitem[MSb]{Monod-Shalom2}
Nicolas Monod and Yehuda Shalom, {\em Orbit equivalence rigidity and bounded
  cohomology}, Preprint.

\bibitem[Nos91]{Noskov}
Gennady~A. Noskov, {\em {Bounded cohomology of discrete groups with
  coefficients}}, Leningr. Math. J. {\bf 2} (1991), no.~5, 1067--1084.

\bibitem[Ol'92]{Olshanskii92}
A.~Yu. Ol'shanski\u{\i}, {\em Almost every group is hyperbolic}, Internat. J.
  Algebra Comput. {\bf 2} (1992), no.~1, 1--17.

\bibitem[Pau96]{Paulin96}
Fr{\'e}d{\'e}ric Paulin, {\em Un groupe hyperbolique est d\'etermin\'e par son
  bord}, J. London Math. Soc. (2) {\bf 54} (1996), no.~1, 50--74.

\bibitem[Rip82]{Rips82}
Eliyahu Rips, {\em Subgroups of small cancellation groups}, Bull. London Math.
  Soc. {\bf 14} (1982), no.~1, 45--47.

\bibitem[Rud91]{RudinFA}
Walter Rudin, {\em Functional analysis}, , International Series in Pure and
  Applied Mathematics, McGraw-Hill Inc., New York, 1991.

\bibitem[Sha00]{Shalom00}
Yehuda Shalom, {\em Rigidity of commensurators and irreducible lattices},
  Invent. Math. {\bf 141} (2000), no.~1, 1--54.

\bibitem[Sha02]{ShalomQI}
Yehuda Shalom, {\em Harmonic analysis, cohomology, and the large scale geometry
  of amenable groups}, Preprint, 2002.

\bibitem[Zim78]{Zimmer78b}
Robert~J. Zimmer, {\em Amenable ergodic group actions and an application to
  {P}oisson boundaries of random walks}, J. Functional Analysis {\bf 27}
  (1978), no.~3, 350--372.

\bibitem[Zim83]{Zimmer83}
Robert~J. Zimmer, {\em Ergodic actions of semisimple groups and product
  relations}, Ann. of Math. (2) {\bf 118} (1983), no.~1, 9--19.

\bibitem[{\.Z}uk02]{ZukERN}
Andrzej {\.Z}uk, {\em On property ({T}) for discrete groups}, Rigidity in
  dynamics and geometry (Cambridge, 2000) (Berlin), Springer, Berlin, 2002,
  pp.~473--482.

\end{thebibliography}
\end{document}